\documentclass[10pt,reqno]{amsart}
\usepackage{amsmath}
\usepackage{amsfonts}
\usepackage{amssymb}
\usepackage{graphicx}
\usepackage{amsthm,graphicx,color,yfonts}
\usepackage{pdfsync}
\usepackage{epstopdf}%
\usepackage{pifont}
\usepackage{bbm}
\usepackage{a4wide}

\usepackage[colorlinks=true]{hyperref}
\hypersetup{linkcolor=red,citecolor=blue,urlcolor=blue} 

\usepackage{wrapfig}
\usepackage{tikz}
\usetikzlibrary{arrows,calc,decorations.pathreplacing}
\definecolor{light-gray1}{gray}{0.90}
\definecolor{light-gray2}{gray}{0.80}
\definecolor{light-gray3}{gray}{0.60}

\newcommand{\R} {\mathbb R}
\newcommand{\cuad}{{\sqcap\kern-.68em\sqcup}}

\newcommand{\be}{\begin{equation}}
\newcommand{\ee}{\end{equation}}

\definecolor{darkgreen}{rgb}{0.2,0.7,0.1}

\renewcommand{\Re}{\mathrm{Re}}
\renewcommand{\Im}{\mathrm{Im}}

\newcommand{\Z}{\mathbb{Z}}
\newcommand{\T}{\mathbb{T}}

\def\bm{\left( \begin{array}{cc}}
\def\endm{\end{array}\right)}

\providecommand{\bra}[1]{\langle #1 \rangle}

\newcommand{\ba}{\begin{equation*}}
\newcommand{\ea}{\begin{equation*}}
\newcommand{\bea}{\begin{eqnarray}}
\newcommand{\eea}{\end{eqnarray}}
\newcommand{\bee}{\begin{eqnarray*}}
\newcommand{\eee}{\end{eqnarray*}}
\newcommand{\ben}{\begin{enumerate}}
\newcommand{\een}{\end{enumerate}}

\numberwithin{equation}{section}

\newtheorem{thm}{Theorem}[section]
\newtheorem*{thm*}{Theorem}
\newtheorem{prop}{Proposition}[section]

\newtheorem{lem}{Lemma}[section]
\newtheorem{definition}{Definition}[section]
\newtheorem{rem}{Remark}[section]

\author{Jaeseop Ahn}
\address{Human-Centered Artificial Intelligence Research Institute, Ewha Womans University, Seoul 03760, Republic of Korea}
\email{j.ahn@ewha.ac.kr}

\author{Chulkwang Kwak}
\address{Department of Mathematics, Ewha Womans University, Seoul 03760, Republic of Korea}
\address{Korea Institute for Advanced Study, Seoul 02455, Republic of Korea}
\email{chkwak@ewha.ac.kr}

\title[Normal form smoothing and spatial analyticity]{Normal form smoothing and spatial analyticity for the dispersion-generalized periodic KdV equation}



\begin{document}

\begin{abstract}
This paper investigates the persistence of spatial analyticity for solutions to the dispersion-generalized periodic Korteweg-de Vries equation
\[\partial_tu - \partial_x|D_x|^\alpha u + u\partial_xu = 0, \qquad \alpha\ge2.\]
For real-valued initial data $u_0 \in G^{\sigma_0,s}(\T)$, with $\sigma_0 > 0$ and $s \ge -\alpha/4$, we prove that the corresponding global solution satisfies
\[u \in C([-T,T];G^{\sigma(T),s}(\T)), \qquad \sigma(T)=\min\left\{\sigma_0,\frac{c}{T}\right\},\]
for every $T>0$. As a consequence, our result improves the lower bound for the radius of spatial analyticity established by Himonas, Kalisch, and Selberg \cite{HKS} for $2\le\alpha<4$, including the classical periodic KdV equation, and establishes a unified $T^{-1}$ lower bound for the analyticity radius for every $\alpha\ge2$.

\medskip

The main novelties are twofold. First, we introduce a normal form reduction at the level of the Gevrey energy, which reveals an additional cancellation in the cubic Gevrey commutator and separates it into a non-accumulating boundary correction and a quartic remainder. Second, we establish a global-in-time periodic $L^4_{t,x}$ estimate for the non-polynomial phase $k|k|^\alpha$, valid for every real $\alpha\ge2$. Together, these ingredients yield a refined almost conservation law and the $T^{-1}$ lower bound for the radius of spatial analyticity.
\end{abstract}

\maketitle

\section{Introduction}\label{sec1}

\subsection{Problem setting}
The Korteweg--de Vries (KdV) equation is one of the fundamental models in the theory of nonlinear dispersive equations. Its origin goes back to Russell's experimental observation of a solitary water wave \cite{Russell}. Boussinesq subsequently developed a mathematical theory for long waves in a shallow channel, identifying the balance between weak nonlinearity and weak dispersion that gives rise to solitary waves \cite{Boussinesq}. Korteweg and de Vries later derived the equation now bearing their names as a model for the propagation of long, weakly nonlinear, and weakly dispersive waves in a shallow channel \cite{KdV}. After a suitable normalization, the classical KdV equation takes the form
\[\partial_tu+\partial_x^3u+u\partial_xu=0.\]
The KdV equation also plays a central role in the theory of completely integrable systems. Gardner, Greene, Kruskal, and Miura introduced the inverse scattering method for the KdV initial-value problem \cite{GGKM}, while Lax formulated its integrability through an isospectral deformation of a Schr\"odinger operator \cite{Lax}. Moreover, the KdV equation possesses a Hamiltonian formulation and an infinite hierarchy of conservation laws \cite{ZF}. For historical accounts of the origin and development of the KdV equation, see, for instance, \cite{Miles, deJager, Hereman} and references therein.

\medskip

In this paper, we consider the Cauchy problem for the dispersion-generalized KdV equation posed on the periodic domain $\T=\R/(2\pi\Z)$:
\begin{equation}\label{eq:dgkdv}
\begin{cases}
\partial_tu - \partial_x|D_x|^\alpha u+u\partial_xu=0,\\
u(0,x)=u_0(x),
\end{cases}
\end{equation}
where $u:\R\times\T\to\R$ and $\alpha\ge2$. Here $D_x=-i\partial_x$ is the Fourier multiplier with symbol $k$, while $|D_x|^\alpha$ is the Fourier multiplier defined by
\[\widehat{|D_x|^\alpha f}(k)=|k|^\alpha\widehat f(k).\]
The corresponding dispersion relation is
\[\omega(k):=k|k|^\alpha.\]
When $\alpha=2$, we have $-\partial_x|D_x|^2=\partial_x^3$, and hence \eqref{eq:dgkdv} reduces to the classical periodic KdV equation.

\medskip

Although the complete integrability of the classical KdV equation does not persist for a general value of $\alpha$, equation \eqref{eq:dgkdv} retains a Hamiltonian structure. Indeed, define
\[\mathcal H_\alpha(u):=\frac12\int_{\T}\big||D_x|^{\alpha/2}u(x)\big|^2\; dx-\frac16\int_{\T}u(x)^3\; dx.\]
Then
\[\frac{\delta\mathcal H_\alpha}{\delta u}=|D_x|^\alpha u-\frac12u^2,\]
and equation \eqref{eq:dgkdv} can formally be written as
\[\partial_tu=\partial_x\frac{\delta\mathcal H_\alpha}{\delta u}.\]
Thus, for sufficiently smooth real-valued solutions, the Hamiltonian $\mathcal H_\alpha(u)$ is formally conserved. The spatial mean and the $L^2$ norm are also conserved. Our argument, however, does not use the complete integrability of the classical KdV equation and applies uniformly to every real $\alpha\ge2$.

\medskip

The notion of well-posedness was introduced by Hadamard \cite{Hadamard1902}. We recall the definition in the Sobolev setting.
\begin{definition}[Local well-posedness]\label{def:well posedness}
Let $s\in\R$ and $u_0\in H^s(\T)$. We say that the Cauchy problem \eqref{eq:dgkdv} is locally well-posed in $H^s(\T)$ if the following conditions are satisfied:
\begin{enumerate}
\item \emph{\bf Existence.} There exist a time $T=T(\|u_0\|_{H^s(\T)})>0$ and a solution space $X_T^s\subset C([-T,T];H^s(\T))$ such that a solution $u$ to \eqref{eq:dgkdv} exists in $X_T^s$.

\medskip

\item \emph{\bf Uniqueness.} The solution is unique in $X_T^s$.

\medskip

\item \emph{\bf Continuous dependence on the initial data.} There exists a neighborhood $B\subset H^s(\T)$ of $u_0$ such that the data-to-solution map
\[\Phi:B\to X_T^s,\qquad \Phi(v_0)=v,\]
is continuous, where $v$ denotes the solution to \eqref{eq:dgkdv} with initial datum $v_0$.
\end{enumerate}
\end{definition}

\begin{rem}
The Cauchy problem is said to be globally well-posed in $H^s(\T)$ if the solution extends uniquely to every finite time interval and depends continuously on the initial data on each such interval.
\end{rem}
We also refer to \cite{Tao2006} for a general discussion of well-posedness for nonlinear dispersive equations.

\medskip

The low-regularity well-posedness theory for the periodic KdV equation was initiated by Bourgain \cite{Bourgain1993-2}, who introduced Fourier restriction spaces and established global well-posedness in $L^2(\T)$. Kenig, Ponce, and Vega \cite{KPV1996} further improved the fundamental bilinear estimate underlying the low-regularity theory. Colliander, Keel, Staffilani, Takaoka, and Tao \cite{CKSTT} proved the sharp global well-posedness in $H^{-1/2}(\T)$. Exploiting the complete integrability of the equation, Kappeler and Topalov \cite{KT} further extended global well-posedness to $H^{-1}(\T)$.

\medskip

For periodic higher-order KdV equations with dispersion of order $2j+1$, Hirayama \cite{Hirayama} established local well-posedness in $H^s(\T)$ for $s\ge-j/2$. In the fifth-order case $j=2$, Kato \cite{KatoKawahara} proved local well-posedness for $s\ge-3/2$ and global well-posedness for $s\ge-1$. For $j\ge3$, Hong and the second author \cite{HK} established global well-posedness in $H^s(\T)$ for $s\ge-j/2$. Li and Shi \cite{LiShi} extended the local well-posedness theory to general real dispersion exponents $\alpha\ge2$ in the range $s\ge-\alpha/4$.

\medskip

We are concerned with the persistence of spatial analyticity for solutions to \eqref{eq:dgkdv}. For $\sigma>0$, define the complex strip
\[S_\sigma:=\{x+iy:x\in\T,\ |y|<\sigma\}.\]
A periodic function $f$ is said to have a uniform radius of spatial analyticity at least $\sigma$ if it admits a periodic holomorphic extension $F$ to $S_\sigma$. When an additional Sobolev regularity index $s\in\R$ is prescribed, we require that
\[\sup_{|y|<\sigma}\|F(\cdot+iy)\|_{H^s(\T)}<\infty.\]
By the periodic Paley-Wiener characterization \cite[Theorem 3.3]{Danielski},\footnote{The result in \cite[Theorem 3.3]{Danielski} is stated for $s\ge0$, while the same conclusion for $s\in\R$ follows by a minor modification of the proof.} this holomorphic extension property together with the uniform $H^s$ bound is equivalent to
\[\sum_{k\in\Z}e^{2\sigma|k|}\bra{k}^{2s}|\widehat f(k)|^2<\infty.\]
For related periodic Fourier characterizations of analyticity, see also \cite[Theorem 2.1]{Rauch} and \cite[Section I.4, Exercise 4]{K}. The corresponding Paley-Wiener theorem on the real line is given in \cite[Theorem 7.2]{K}. This characterization leads naturally to the Gevrey space $G^{\sigma,s}(\T)$, whose precise definition will be given in Section \ref{sec:prelim}.

\medskip

Suppose that the initial datum $u_0$ has a uniform radius of analyticity $\sigma_0>0$. The local Gevrey result stated in \cite[Theorem 4]{HKS} ensures that the corresponding solution remains analytic for a short time with the same radius. The global persistence problem asks whether there exists a positive function $\sigma(t)$ such that
\[u(t)\in G^{\sigma(t),s}(\T)\]
for every later time. A more quantitative problem is to determine how rapidly the available radius $\sigma(t)$ may decrease as $|t|\to\infty$.

\medskip

This type of problem was first studied in an abstract framework by Kato and Masuda \cite{KM} and has since received considerable attention for nonlinear dispersive equations. For the KdV equation on the real line, Bona and Gruji\'c \cite{BG} showed that the radius of analyticity cannot decay faster than $e^{-ct^2}$. Bona, Gruji\'c, and Kalisch \cite{BGK} subsequently obtained the algebraic lower bound
\[\sigma(t)\gtrsim |t|^{-12}.\]
Selberg and da Silva \cite{SS} further improved this result and proved that, for every $\epsilon>0$,
\[\sigma(t)\gtrsim_\epsilon |t|^{-4/3-\epsilon}.\]
Tesfahun \cite{Te2} removed the $\epsilon$-loss and obtained the lower bound of order $|t|^{-4/3}$. The rate was further improved by Huang and Wang \cite{HW} to
\[\sigma(t)\gtrsim |t|^{-1/4}.\]

\medskip

Related results have been obtained for the generalized KdV equation \cite{ST_gKdV, BP}, dissipative and weakly damped KdV equations \cite{LW, W}, and higher-order nonlinear dispersive equations \cite{ZLD}. For other nonlinear dispersive equations and systems, see, for instance, \cite{BGK2, ST, Te, PS, S, AKS1, AKS2, FH, AKS3, TTB, BG_NLS, FNP, FP1, FP2, BC, Sato} and the references therein.

\medskip

Compared with the corresponding problems on the real line, considerably less is known about quantitative persistence of spatial analyticity on periodic domains. For the periodic KdV equation, Kappeler and P\"oschel \cite{KP} used global Birkhoff coordinates to show that analytic solutions remain uniformly bounded for all time in a fixed, possibly smaller, analytic strip, yielding a time independent positive lower bound for the radius of spatial analyticity. Their argument relies essentially on the complete integrability of the equation. A related result for the defocusing cubic nonlinear Schr\"odinger equation on the circle was obtained by Maspero \cite{Maspero}.

\medskip

Beyond these integrable settings, however, it remains natural to seek general analytic methods, such as harmonic analysis and functional analysis, for quantitative persistence of spatial analyticity. For the periodic dispersion-generalized KdV equation \eqref{eq:dgkdv}, Himonas, Kalisch, and Selberg \cite{HKS} established the asymptotic lower bound
\[\sigma(t)\ge c|t|^{-p},\qquad p=\max\left\{1,\frac4\alpha\right\},\]
as $|t|\to\infty$. In particular, their result gives a lower bound of order $|t|^{-2}$ for the classical periodic KdV equation. For the periodic Benjamin--Bona--Mahony equation, Himonas and Petronilho \cite{HP} obtained a lower bound of order $|t|^{-1}$, which was subsequently improved by Getachew \cite{Get1} to $|t|^{-2/3}$. Getachew \cite{Get2} also obtained the same lower bound for the periodic fifth-order KdV--BBM equation.

\subsection{Main result}
Motivated by these works, we establish a uniform algebraic lower bound of order $|t|^{-1}$ for the radius of spatial analyticity of solutions to \eqref{eq:dgkdv} for every $\alpha\ge2$. Our main result is stated as follows.

\begin{thm}\label{thm:main}
Let $\alpha\ge2$, $\sigma_0>0$, and $s\ge-\alpha/4$, and suppose that the initial datum $u_0\in G^{\sigma_0,s}(\T)$ is real-valued. Then there exists a global solution $u$ to \eqref{eq:dgkdv} such that, for every $T>0$,
\[u\in C([-T,T];G^{\sigma(T),s}(\T)),\]
where
\[\sigma(T):=\min\left\{\sigma_0,\frac cT\right\}\]
for some constant $c>0$ depending on $\|u_0\|_{G^{\sigma_0,s}}$, $\alpha$, $\sigma_0$, and $s$. The solution is unconditionally unique in $C([-T,T];G^{\sigma(T),s}(\T))$ in the sense of Remark \ref{rem:unconditional uniqueness}. Moreover, the solution depends continuously on the initial datum on every finite time interval in the sense of Remark \ref{rem:continuous dependence}.
\end{thm}

\begin{rem}
For $2\le\alpha<4$, Theorem \ref{thm:main} improves the lower bound obtained in \cite{HKS} from order $T^{-4/\alpha}$ to order $T^{-1}$. In particular, for the classical periodic KdV equation $\alpha=2$, the lower bound of order $T^{-2}$ obtained in \cite{HKS} is improved to order $T^{-1}$. For $\alpha\ge4$, the bound in \cite{HKS} is already of order $T^{-1}$, and Theorem \ref{thm:main} recovers the same rate through a unified normal form argument valid for every $\alpha\ge2$.
\end{rem}

\begin{rem}
The restriction $s\ge-\alpha/4$ comes from the sharp $X^{s,b}$ bilinear estimate used in the local well-posedness argument. Nevertheless, global existence and persistence extend to every $s\in\R$ after a fixed reduction of the initial analytic radius.

More precisely, let $\sigma_0>0$, $s\in\R$, and let $u_0\in G^{\sigma_0,s}(\T)$ be real-valued. Fix $\sigma_\ast\in(0,\sigma_0)$ and choose $\rho_0$ such that
\[\sigma_\ast<\rho_0<\sigma_0.\]
The loss of analytic radius gives
\[G^{\sigma_0,s}(\T)\hookrightarrow G^{\rho_0,0}(\T).\]
Applying the case $s=0$ of Theorem \ref{thm:main} to the same initial datum $u_0$, regarded as an element of $G^{\rho_0,0}(\T)$, yields a global solution $u$ to \eqref{eq:dgkdv} satisfying
\[u\in C([-T,T];G^{\rho(T),0}(\T))\]
for every $T>0$, where
\[\rho(T):=\min\left\{\rho_0,\frac{c_{\rho_0}}{T}\right\}\]
for some constant $c_{\rho_0}>0$ depending on $\|u_0\|_{G^{\sigma_0,s}}$, $\alpha$, $\sigma_0$, $s$, and $\rho_0$.

Choose $c_\ast\in(0,c_{\rho_0})$ and set
\[\sigma_\ast(T):=\min\left\{\sigma_\ast,\frac{c_\ast}{T}\right\}.\]
Then $\sigma_\ast(T)<\rho(T)$ for every $T>0$. Hence,
\[G^{\rho(T),0}(\T)\hookrightarrow G^{\sigma_\ast(T),s}(\T),\]
and therefore
\[u\in C([-T,T];G^{\sigma_\ast(T),s}(\T))\]
for every $T>0$. Thus, global existence and persistence hold for every $s\in\R$ after a fixed loss of analytic radius. Moreover, the uniform radius of spatial analyticity satisfies
\[\sigma(t)\gtrsim |t|^{-1}\]
as $|t|\to\infty$. This does not assert well-posedness in $G^{\sigma_0,s}(\T)$ when $s<-\alpha/4$.
\end{rem}

The following local result from \cite{HKS} is the local input for the proof.

\begin{thm}[Theorem 4 in \cite{HKS}]\label{thm:hks lwp}
Let $\alpha\ge2$, $\sigma>0$, and $s\ge-\alpha/4$. For every real-valued initial datum $u_0\in G^{\sigma,s}(\T)$, there exists
\[\delta=c_0\left(1+\|u_0\|_{G^{\sigma,s}}\right)^{-a}>0,\]
where $a,c_0>0$ depend only on $\alpha$ and $s$, such that \eqref{eq:dgkdv} has a unique solution
\[u\in C([-\delta,\delta];G^{\sigma,s}(\T)).\]
Moreover,
\[\|u\|_{Y_{[-\delta,\delta]}^{\sigma,s}}\le C\|u_0\|_{G^{\sigma,s}},\]
where $C>0$ depends only on $\alpha$ and $s$.
\end{thm}

The space $Y_I^{\sigma,s}$ will be defined in Section \ref{sec:prelim}. A crucial feature of Theorem \ref{thm:hks lwp} is that the constants in the lifespan and the local solution estimate are independent of the analytic radius $\sigma$. This uniformity allows us to iterate the local theory while decreasing the analytic radius. We note that the proof in \cite{HKS} is carried out after reducing to the zero-mean case by a Galilean transformation. We therefore use the $Y^{\sigma,s}$ estimate below for mean-zero solutions, and treat the non-zero mean case separately through the corresponding shifted dispersion relation.

\begin{rem}\label{rem:unconditional uniqueness}
We justify the unconditional uniqueness assertion in Theorem \ref{thm:main}. In fact, the argument applies to every positive-radius Gevrey class. Here, a distributional solution in $C(I;G^{\sigma,s}(\T))$ means a function in this space that satisfies \eqref{eq:dgkdv} in $\mathcal D'(I^\circ\times\T)$. Let $I\subset\R$ be a bounded interval, and suppose that
\[u,v\in C(I;G^{\sigma,s}(\T))\]
are two distributional solutions to \eqref{eq:dgkdv} satisfying $u(t_0)=v(t_0)$ for some $t_0\in I$. Setting $w:=u-v$, we obtain
\[\partial_tw-\partial_x|D_x|^\alpha w+\frac12\partial_x((u+v)w)=0.\]
Since $\sigma>0$, we have
\[G^{\sigma,s}(\T)\hookrightarrow H^m(\T)\]
for every $m\in\R$. Choosing $m$ sufficiently large, we obtain
\[\partial_x|D_x|^\alpha w\in C(I;L^2(\T)),\qquad \partial_x((u+v)w)\in C(I;L^2(\T)),\qquad \partial_x(u+v)\in C(I;L^\infty(\T)).\]
The difference equation therefore implies $\partial_tw\in C(I;L^2(\T))$, and hence $w\in C^1(I;L^2(\T))$. Since $\partial_x|D_x|^\alpha$ is skew-adjoint, we obtain
\[\frac{d}{dt}\|w(t)\|_{L^2}^2
=-\frac12\int_{\T}\partial_x(u+v)(t,x)|w(t,x)|^2\,dx.\]
Consequently,
\[\left|\frac{d}{dt}\|w(t)\|_{L^2}^2\right|
\le\frac12\|\partial_x(u+v)(t)\|_{L^\infty}\|w(t)\|_{L^2}^2.\]
Since $w(t_0)=0$, Gr\"onwall's inequality applied on both sides of $t_0$ implies
\[\|w(t)\|_{L^2}=0\]
for every $t\in I$. By Plancherel's theorem,
\[\widehat w(t,k)=0\]
for every $k\in\Z$ and $t\in I$. It follows that
\[\|w(t)\|_{G^{\sigma,s}}=0\]
for every $t\in I$. Thus, $u=v$ in $C(I;G^{\sigma,s}(\T))$.
\end{rem}

\begin{rem}\label{rem:continuous dependence}
Although continuous dependence is not stated explicitly in Theorem \ref{thm:hks lwp}, the contraction argument used in its proof, together with the Galilean transformation, yields local Lipschitz continuity on bounded sets of initial data with a common spatial mean. 

\medskip

More precisely, for every $A>0$, there exists $\delta_A>0$, depending only on $\alpha$, $s$, and $A$, such that $\delta_A$ and the implicit constant below are independent of $\sigma$. If $u_0,v_0\in G^{\sigma,s}(\T)$ are real-valued and satisfy
\[\|u_0\|_{G^{\sigma,s}}\le A,\qquad \|v_0\|_{G^{\sigma,s}}\le A,\qquad \widehat u_0(0)=\widehat v_0(0),\]
then the corresponding solutions satisfy
\[\|u-v\|_{C([-\delta_A,\delta_A];G^{\sigma,s}(\T))}\lesssim_{\alpha,s,A}\|u_0-v_0\|_{G^{\sigma,s}}.\]
Indeed, the same Galilean transformation reduces both solutions to mean-zero solutions, for which the contraction argument gives the corresponding estimate in $Y_{[-\delta_A,\delta_A]}^{\sigma,s}$, and spatial translations are isometries on $G^{\sigma,s}(\T)$.

\medskip

For general initial data, let $f_n\to f$ in $G^{\sigma,s}(\T)$. Since the sequence is bounded, fix $A>0$ such that
\[\|f_n\|_{G^{\sigma,s}}\le A,\qquad \|f\|_{G^{\sigma,s}}\le A\]
for every $n$, and set
\[m_n:=\widehat f_n(0),\qquad m:=\widehat f(0),\qquad f_n^\circ:=f_n-m_n,\qquad f^\circ:=f-m.\]
Then
\[m_n\to m,\qquad f_n^\circ\to f^\circ\quad\text{in}\quad G^{\sigma,s}(\T).\]
Let $v_n$ and $v$ denote the corresponding mean-zero solutions. The original solutions are given by
\[u_n(t,x)=v_n(t,x-m_nt)+m_n,\qquad u(t,x)=v(t,x-mt)+m.\]
Writing $\tau_yg(x):=g(x-y)$, we have
\[\sup_{|t|\le\delta_A}\|u_n(t)-u(t)\|_{G^{\sigma,s}} \le \sup_{|t|\le\delta_A}\|v_n(t)-v(t)\|_{G^{\sigma,s}} +\sup_{|t|\le\delta_A}\|\tau_{(m_n-m)t}v(t)-v(t)\|_{G^{\sigma,s}} +|m_n-m|.\]
The first term tends to zero by the local Lipschitz estimate for mean-zero data. For the second term, since spatial translations form a strongly continuous group of isometries on $G^{\sigma,s}(\T)$, the map
\[\Psi:\R\times G^{\sigma,s}(\T)\to G^{\sigma,s}(\T),\qquad \Psi(y,g):=\tau_yg,\]
is jointly continuous. Since
\[K:=\{v(t):|t|\le\delta_A\}\]
is compact, the Heine--Cantor theorem implies that $\Psi$ is uniformly continuous on $[-1,1]\times K$. Comparing $\Psi(y,g)$ with $\Psi(0,g)$ uniformly for $g\in K$, we obtain
\[\sup_{g\in K}\|\tau_yg-g\|_{G^{\sigma,s}} \to 0\]
as $y\to0$. Since
\[\sup_{|t|\le\delta_A}|(m_n-m)t|\le\delta_A|m_n-m|\to0,\]
the second term also tends to zero. The last term tends to zero because $m_n\to m$. Hence,
\[u_n\to u\quad\text{in}\quad C([-\delta_A,\delta_A];G^{\sigma,s}(\T)).\]

\medskip

Continuous dependence on a finite time interval is understood with an analytic radius chosen uniformly on bounded sets of initial data. Given $A>0$, define
\[B_A^{\sigma_0,s}:=\left\{f\in G^{\sigma_0,s}(\T):f\text{ is real-valued and }\|f\|_{G^{\sigma_0,s}}\le A\right\}.\]
For every $T>0$, the local lifespan, the local solution bounds, and the lower bound for the analytic radius used in the proof of Theorem \ref{thm:main} can be chosen uniformly for initial data in $B_A^{\sigma_0,s}$. Consequently, there exists $c_A>0$, depending only on $A$, $\alpha$, $\sigma_0$, and $s$, such that, with
\[\sigma_A(T):=\min\left\{\sigma_0,\frac{c_A}{T}\right\},\]
the data-to-solution map
\[\Phi_T:B_A^{\sigma_0,s}\to C([-T,T];G^{\sigma_A(T),s}(\T)),\qquad \Phi_T(f)=u_f,\]
is continuous. Indeed, the uniform bounds obtained in the proof of Theorem \ref{thm:main} keep the intermediate data in a common bounded subset of $G^{\sigma_A(T),s}(\T)$. One may therefore iterate the preceding local continuous-dependence argument over the finitely many local intervals used in the proof of Theorem \ref{thm:main}.
\end{rem}

A central approach to obtaining quantitative lower bounds for the radius of spatial analyticity is to identify a nonlinear smoothing effect in the evolution of the Gevrey energy. On the real line, the required gain can be obtained through frequency-localized dispersive estimates, including local smoothing and maximal-function estimates, see for instance, \cite{Te2, BC}, and references therein. The corresponding real-line smoothing mechanisms are not available on the periodic domain. To overcome this difficulty, Himonas, Kalisch, and Selberg \cite{HKS} controlled the Gevrey commutator at the $G^{\sigma,0}$ energy level by means of $X^{s,b}$ bilinear estimates at negative Sobolev regularity, thereby exploiting the regularity gap between the energy and multilinear estimates.

One of the main novelties of this paper is the use of a normal form reduction to extract the required smoothing directly from the algebraic structure of the interaction phase. More precisely, we apply the normal form reduction to the evolution of the Gevrey energy, revealing an alternative nonlinear smoothing mechanism encoded in the nonresonant interaction phase.

The normal form method goes back to Shatah's work on quadratic nonlinear Klein--Gordon equations \cite{Shatah1985}. In the periodic KdV setting, Babin, Ilyin, and Titi used a successive averaging procedure to reveal the regularization mechanism of the equation \cite{BIT2011}. Differentiation-by-parts and normal form techniques have since been applied to various periodic dispersive equations, see, for instance, \cite{KO, GKO, KwakKawahara, KY}.

We now describe the normal form reduction used in this paper. We introduce the interaction representation
\[v(t):=e^{-it\omega(D_x)}u(t)\]
and the Gevrey profile
\[V(t):=e^{\sigma|D_x|}v(t).\]
The corresponding Gevrey energy is
\[E_\sigma(t):=\sum_{k\in\Z}|\widehat V(t,k)|^2=\|u(t)\|_{G^{\sigma,0}}^2.\]
Differentiating this energy produces a cubic oscillatory expression associated with the cubic phase
\[\Phi_3(k_1,k_2,k_3):=\omega(k_1)+\omega(k_2)+\omega(k_3),\qquad k_1+k_2+k_3=0.\]

A central algebraic ingredient is the factorization
\[\Phi_3(k_1,k_2,k_3)=k_1k_2k_3Q_3(k_1,k_2,k_3),\]
where
\[|Q_3(k_1,k_2,k_3)|\sim_\alpha K_{\max}^{\alpha-2},\qquad K_{\max}:=\max\{|k_1|,|k_2|,|k_3|\}.\]
Under the mean-zero assumption, interactions involving the zero frequency are absent. Consequently, there are no non-trivial cubic resonances, and $\Phi_3$ is nonzero on the cubic interaction set. The latter property allows us to write the oscillatory factor as
\[e^{it\Phi_3}=-\frac i{\Phi_3}\partial_t(e^{it\Phi_3}).\]
An integration by parts in time transforms the cubic term in the Gevrey energy identity into a cubic boundary term $B(t)$ and a quartic remainder $R_4(t)$. This is the normal form reduction carried out in Section \ref{sec:normal form}.

The phase denominator provides the smoothing needed in the estimates for both terms generated by the normal form reduction. The cubic boundary term satisfies
\[|B(t)|\lesssim_\alpha\sigma E_\sigma(t)^{3/2}.\]
Since it occurs only at the endpoints of the full time interval, it does not accumulate at each step of the local iteration. On the other hand, if $\delta\in(0,1]$ denotes the local lifespan, then the quartic remainder satisfies
\begin{equation}\label{eq:quartic estimate}
\left|\int_0^tR_4(t')\; dt'\right|\lesssim_\alpha\sigma\delta^{1/4}\|u\|_{Y_\delta^{\sigma,0}}^4,\qquad 0\le t\le\delta.
\end{equation}
This term accumulates over successive local intervals, leading to the condition $T\sigma\lesssim1$ and hence to the lower bound of order $T^{-1}$.

\begin{rem}
The normal form reduction used here is closely related to the higher modified energy construction used in the $I$-method. For the KdV equation on the real line, Huang and Wang \cite{HW} introduced successive correction terms to the Gevrey energy, with the corresponding multipliers obtained by dividing the preceding energy increments by the appropriate resonance functions. At the level of the first correction, this modified energy construction may be viewed as the energy-level counterpart of an integration by parts in time in the normal form method. The distinction at this stage is mainly one of formulation: the $I$-method encodes the cancellation through a corrected energy, whereas our argument exhibits it directly in the oscillatory Gevrey energy identity.
\end{rem}

A second novelty of this paper is the derivation of a global-in-time periodic $L^4$ estimate for arbitrary real dispersion exponents $\alpha\ge2$. Periodic $L^4$ estimates play a fundamental role in the low-regularity well-posedness theory for dispersive equations on the torus. For representative applications, see, for instance, \cite{Bourgain1993-1, Bourgain1993-2, Tao2001, TakaokaTsutsumi2004, Kwak2018, CKV2022}. The periodic estimates in these works concern integer-order dispersive phases, whereas our result also applies to the non-polynomial phase $\omega(k)=k|k|^\alpha$ when $\alpha$ is non-integer. Note that noninteger dispersion relations in the periodic setting have also been considered, in a different context, by Demirbas, Erdo\u{g}an, and Tzirakis \cite{DET}. In the present setting, we prove that
\[\|u\|_{L^4_{t,x}(\R\times\T)}\lesssim_{\alpha,b}\|u\|_{X^{0,b}}\]
whenever
\[b>\frac{\alpha+2}{4(\alpha+1)}.\]
Since this threshold is at most $1/3$ for every $\alpha\ge2$, smooth time localization yields the factor $\delta^{1/4}$ in \eqref{eq:quartic estimate}. This gain is essential for controlling the accumulation of the quartic remainder over the local iteration.

\medskip

The normal form identity therefore takes the form
\[E_\sigma(t)=E_\sigma(0)+B(t)-B(0)+\int_0^tR_4(t')\; dt'.\]
The preceding estimates, together with the local theory, close a continuity argument for $\sigma$ of order $T^{-1}$ and prove the desired lower bound for $s=0$.

The reduction to the case $s=0$ is standard, see \cite{SS}. For completeness, we include the short argument. Indeed, sacrificing a fixed portion of the analytic radius gives
\[\|f\|_{G^{\rho/2,s}}\lesssim_{\rho,s}\|f\|_{G^{\rho,0}}.\]

The mean-zero assumption used in the normal form argument causes no loss of generality. The general case is recovered by the Galilean transformation in Remark \ref{rem:recover mean}.

\begin{rem}
The power of $\sigma$ gained in an almost conservation law determines the resulting lower bound for the radius of analyticity. The exponential weight $e^{\sigma|D_x|}$ yields at most one power of $\sigma$ through the direct commutator estimate. In \cite{HKS}, the amount of this gain that can be used is further restricted by the available negative-regularity bilinear estimate. As the dispersion becomes stronger, this bilinear estimate provides a larger regularity gain and eventually compensates for the entire derivative loss.

In the present paper, the normal form reduction provides the required smoothing directly through the phase denominator. This compensates for the derivative loss for every $\alpha\ge2$ and allows us to retain the full factor $\sigma$. The resulting almost conservation law therefore leads to the lower bound of order $T^{-1}$ throughout the entire range $\alpha\ge2$.
\end{rem}

\begin{rem}
The lower bound of order $T^{-1}$ reflects the limitation of the exponential weight in the present argument rather than an intrinsic limitation of the Gevrey class. Even analytic weights such as $\cosh(\sigma|D_x|)$, or a more detailed use of the Taylor expansion of the analytic weight, may yield more than one power of $\sigma$, see, for instance, \cite{HW, FP2026}. The additional gain in $\sigma$, however, is accompanied by an additional derivative loss. On the real line, this loss may be handled using local smoothing and maximal-function estimates, whereas the corresponding real-line smoothing mechanism is unavailable on the periodic domain.

For $\alpha>2$, if the exponential weight is replaced by the hyperbolic cosine weight $\cosh(\sigma|D_x|)$, the normal form denominator provides a surplus smoothing of order $\alpha-2$. It is therefore not necessary to extract the full quadratic gain in $\sigma$ available from this weight. One may instead expect to trade only the available surplus smoothing for an additional fractional power of $\sigma$. This suggests an almost conservation law with gain $\sigma^{1+\min\{1,\alpha-2\}}$, and hence a lower bound of order $T^{-1/(\alpha-1)}$ for $2<\alpha<3$ and of order $T^{-1/2}$ for $\alpha\ge3$. We do not pursue these refinements in the present paper.
\end{rem}

\begin{rem}[Application to the periodic modified KdV equation]
The normal form approach developed in the present paper also extends, in a substantially different form, to the periodic modified KdV equation. In ongoing work \cite{AK2026}, we exploit the exact cancellation associated with the cubic resonance structure and apply a normal form reduction to the remaining nonresonant interactions. This leads to an additional smoothing effect, which compensates for the derivative loss in the Gevrey commutator and yields a refined approximate conservation law. As a consequence, we obtain an improved algebraic lower bound for the radius of spatial analyticity.
\end{rem}

\subsection*{Paper organization}
The remainder of the paper is organized as follows. In Section \ref{sec:prelim}, we introduce the notation and function spaces used throughout the paper and establish a new global-in-time periodic $L^4$ estimate for the non-polynomial dispersion relation. In Section \ref{sec:normal form}, we derive the normal form identity. In Section \ref{sec:estimates}, we establish the required estimates. Finally, in Section \ref{sec:proof}, we prove Theorem \ref{thm:main}.

\subsection*{Acknowledgments}
Both authors were supported by the Global--Learning \& Academic Research Institution for Master's·PhD Students, and Postdocs (G-LAMP) Program of the National Research Foundation of Korea (NRF), funded by the Ministry of Education (No. RS-2025-25442252). C. K. was also partially supported by the Young Research Program of the NRF, funded by the Korean government (MSIT) (No. RS-2023-00210210).

\section{Preliminaries}\label{sec:prelim}

Throughout the paper, we write $A\lesssim_\theta B$ if there exists a constant $C(\theta)>0$ such that
\[A\le C(\theta)B.\]
We write $A\sim_\theta B$ if $A\lesssim_\theta B$ and $B\lesssim_\theta A$. When the dependence of the implicit constant is unimportant, the subscript is omitted. For $y\in\R$, we use the Japanese bracket
\[\langle y\rangle:=(1+|y|^2)^{1/2}.\]

We consider functions on the spatial torus
\[\T:=\R/(2\pi\Z).\]
For an integrable function $f$ on $\T$, define
\[\widehat f(k):=\frac1{2\pi}\int_{\T}e^{-ikx}f(x)\; dx,\qquad k\in\Z.\]
The corresponding Fourier series is
\[f(x)=\sum_{k\in\Z}\widehat f(k)e^{ikx}.\]
With this convention,
\[\widehat{fg}(k)=\sum_{k_1+k_2=k}\widehat f(k_1)\widehat g(k_2)\]
and
\[\|f\|_{L^2(\T)}^2=2\pi\sum_{k\in\Z}|\widehat f(k)|^2.\]

For a space-time function $u$ on $\R\times\T$, define
\[\widetilde u(\tau,k):=\int_{\R}e^{-it\tau}\widehat u(t,k)\; dt=\frac1{2\pi}\int_{\R}\int_{\T}e^{-i(t\tau+xk)}u(t,x)\; dxdt.\]
The corresponding inversion formula is
\[u(t,x)=\frac1{2\pi}\sum_{k\in\Z}\int_{\R}e^{i(t\tau+xk)}\widetilde u(\tau,k)\; d\tau.\]

Set
\[\omega(k):=k|k|^\alpha.\]
We denote by $\omega(D_x)$ the Fourier multiplier with symbol $\omega(k)$.

For real numbers $s$ and $b$, the Fourier restriction space $X^{s,b}=X^{s,b}(\R\times\T)$ associated with the dispersion relation $\omega$ is defined by
\[\|u\|_{X^{s,b}}:=\left(\sum_{k\in\Z}\int_{\R}\langle k\rangle^{2s}\langle\tau-\omega(k)\rangle^{2b}|\widetilde u(\tau,k)|^2\; d\tau\right)^{1/2}.\]
The $X^{s,b}$ space was introduced in its present form by Bourgain \cite{Bourgain1993-1,Bourgain1993-2} in his work on the periodic NLS and generalized KdV equations, although related spaces had already appeared in Beals' work on one-dimensional wave equations \cite{Beals}. Since then, these spaces have played a central role in the analysis of dispersive equations and have been further developed by many authors, including Kenig, Ponce, and Vega \cite{KPV1996} and Tao \cite{Tao2001}.

The periodic KdV iteration requires an endpoint refinement of $X^{s,1/2}$. Colliander, Keel, Staffilani, Takaoka, and Tao \cite{CKSTT} introduced an additional $\ell_k^2L_\tau^1$ component to ensure continuity of the solution in time, while Kenig, Ponce, and Vega \cite{KPV1996} showed that the periodic bilinear estimate relevant to the low-regularity KdV iteration fails unless $b=1/2$. Since $X^{s,1/2}$ alone does not embed into $C(\R;H^s(\T))$, we define
\[\|u\|_{Y^s}:=\|u\|_{X^{s,1/2}}+\left(\sum_{k\in\Z}\langle k\rangle^{2s}\left(\int_{\R}|\widetilde u(\tau,k)|\; d\tau\right)^2\right)^{1/2}.\]
The additional term restores the required continuity in time. By contrast, on the real line one may work with $b>1/2$. Although the stronger dispersive smoothing for $\alpha>2$ may also permit estimates with $b>1/2$ in the periodic setting, we use $Y^s$ uniformly for every $\alpha\ge2$ in order to treat the classical and higher-order cases within the same framework.

For $\sigma\ge0$ and $s\in\R$, define the Gevrey space
\[G^{\sigma,s}(\T):=\left\{f\in\mathcal D'(\T):\|f\|_{G^{\sigma,s}}<\infty\right\},\]
where
\[\|f\|_{G^{\sigma,s}}:=\left(\sum_{k\in\Z}e^{2\sigma|k|}\langle k\rangle^{2s}|\widehat f(k)|^2\right)^{1/2}.\]
When $\sigma=0$, the space $G^{0,s}(\T)$ coincides with the Sobolev space $H^s(\T)$ with an equivalent norm.

For $\sigma\ge0$ and $s,b\in\R$, define the Gevrey-modified Fourier restriction space $X^{\sigma,s,b}$ by
\[\|u\|_{X^{\sigma,s,b}}:=\|e^{\sigma|D_x|}u\|_{X^{s,b}}.\]
Equivalently,
\[\|u\|_{X^{\sigma,s,b}}=\left(\sum_{k\in\Z}\int_{\R}e^{2\sigma|k|}\langle k\rangle^{2s}\langle\tau-\omega(k)\rangle^{2b}|\widetilde u(\tau,k)|^2\; d\tau\right)^{1/2}.\]
In particular,
\[X^{0,s,b}=X^{s,b}.\]

Similarly, define the Gevrey-modified space $Y^{\sigma,s}$ by
\[\|u\|_{Y^{\sigma,s}}:=\|e^{\sigma|D_x|}u\|_{Y^s}.\]
Equivalently,
\[\|u\|_{Y^{\sigma,s}}=\|u\|_{X^{\sigma,s,1/2}}+\left(\sum_{k\in\Z}e^{2\sigma|k|}\langle k\rangle^{2s}\left(\int_{\R}|\widetilde u(\tau,k)|\; d\tau\right)^2\right)^{1/2}.\]
In particular,
\[Y^{0,s}=Y^s.\]

For a time interval $I\subset\R$, define the restriction norms by
\[\|u\|_{X_I^{\sigma,s,b}}:=\inf\left\{\|v\|_{X^{\sigma,s,b}}:v\in X^{\sigma,s,b},\ v=u\text{ on }\operatorname{Int}(I)\times\T\right\}\]
and
\[\|u\|_{Y_I^{\sigma,s}}:=\inf\left\{\|v\|_{Y^{\sigma,s}}:v\in Y^{\sigma,s},\ v=u\text{ on }\operatorname{Int}(I)\times\T\right\}.\]
We use the abbreviations
\[X_I^{s,b}:=X_I^{0,s,b},\qquad Y_I^s:=Y_I^{0,s}.\]
When $I=[0,\delta]$, we write
\[X_\delta^{\sigma,s,b}:=X_{[0,\delta]}^{\sigma,s,b},\qquad Y_\delta^{\sigma,s}:=Y_{[0,\delta]}^{\sigma,s}.\]

\begin{lem}
Let $\sigma\ge0$ and $s\in\R$. Then
\[Y^{\sigma,s}\hookrightarrow C(\R;G^{\sigma,s}(\T))\]
and
\[\sup_{t\in\R}\|u(t)\|_{G^{\sigma,s}}\lesssim\|u\|_{Y^{\sigma,s}}.\]
Consequently, for every time interval $I$,
\[Y_I^{\sigma,s}\hookrightarrow C(I;G^{\sigma,s}(\T)).\]
\end{lem}

\begin{proof}
The proof follows from \cite[Lemma 1]{HKS}. The result for restriction spaces follows by the standard argument of taking the infimum over all extensions.
\end{proof}

\begin{lem}[Smooth time localization]\label{lem:time localization}
Let $\eta\in\mathcal S(\R)$, $\sigma\ge0$, $s\in\R$, $-1/2<b_0\le b_1<1/2$, and $0<T\le1$. Then
\[\|\eta(t/T)u\|_{X^{\sigma,s,b_0}}\lesssim_{\eta,b_0,b_1}T^{b_1-b_0}\|u\|_{X^{\sigma,s,b_1}}.\]
\end{lem}

\begin{proof}
The case $\sigma=0$ follows from \cite[Lemma 2.11]{Tao2006}. Since $e^{\sigma|D_x|}$ commutes with multiplication by $\eta(t/T)$, applying the same estimate to $e^{\sigma|D_x|}u$ gives
\[\|\eta(t/T)u\|_{X^{\sigma,s,b_0}}\lesssim_{\eta,b_0,b_1}T^{b_1-b_0}\|u\|_{X^{\sigma,s,b_1}}.\]
\end{proof}

\begin{lem}[Interval truncation]\label{lem:interval truncation}
Let $\sigma\ge0$, $s\in\R$, and $-1/2<b<1/2$. Then, for every bounded interval $I\subset\R$,
\[\|\mathbf 1_Iu\|_{X^{\sigma,s,b}}\lesssim_b\|u\|_{X^{\sigma,s,b}}.\]
The implicit constant is independent of the position and length of $I$.
\end{lem}

\begin{proof}
We first prove
\begin{equation}\label{eq:interval Hb multiplier}
\|\mathbf 1_If\|_{H^b(\R)}\lesssim_b\|f\|_{H^b(\R)}
\end{equation}
uniformly over all bounded intervals $I\subset\R$.

For a nonnegative locally integrable weight $\varpi$, let $L^2(\varpi)$ denote the weighted $L^2$ space equipped with the norm
\[\|g\|_{L^2(\varpi)}:=\left(\int_{\R}|g(\xi)|^2\varpi(\xi)\; d\xi\right)^{1/2}.\]
Recall that $\varpi$ belongs to the Muckenhoupt class $A_2(\R)$ if
\[[\varpi]_{A_2}:=\sup_{J\subset\R}\left(\frac1{|J|}\int_J\varpi(\xi)\; d\xi\right)\left(\frac1{|J|}\int_J\varpi(\xi)^{-1}\; d\xi\right)<\infty,\]
where the supremum is taken over all bounded intervals $J\subset\R$. By \cite[Theorem 9]{HMW}, this condition is equivalent to the boundedness of the Hilbert transform $\mathcal H$ on $L^2(\varpi)$, where $\mathcal H$ is defined by
\[\mathcal Hg(\xi):=\frac1\pi\operatorname{p.v.}\int_{\R}\frac{g(\eta)}{\xi-\eta}\; d\eta.\]

Set
\[\varpi_b(\xi):=\langle\xi\rangle^{2b}.\]
We claim that $\varpi_b\in A_2(\R)$ whenever $-1/2<b<1/2$. Let $J\subset\R$ be an interval of length $L$. Suppose first that $\operatorname{dist}(J,0)\ge L$. Then
\[\langle\xi\rangle\sim\langle\eta\rangle\]
for all $\xi,\eta\in J$, and hence
\[\left(\frac1{|J|}\int_J\varpi_b(\xi)\; d\xi\right)\left(\frac1{|J|}\int_J\varpi_b(\xi)^{-1}\; d\xi\right)\lesssim_b1.\]

Suppose next that $\operatorname{dist}(J,0)<L$. Then $J\subset[-2L,2L]$. If $L\le1$, both $\varpi_b$ and $\varpi_b^{-1}$ are uniformly bounded on $J$. If $L>1$, the condition $-1<2b<1$ gives
\[\int_J\langle\xi\rangle^{2b}\; d\xi\lesssim_bL^{1+2b},\qquad \int_J\langle\xi\rangle^{-2b}\; d\xi\lesssim_bL^{1-2b}.\]
Thus, in either case,
\[\left(\frac1{|J|}\int_J\varpi_b(\xi)\; d\xi\right)\left(\frac1{|J|}\int_J\varpi_b(\xi)^{-1}\; d\xi\right)\lesssim_b1.\]
Therefore $\varpi_b\in A_2(\R)$.

Let
\[H(t):=\mathbf 1_{(0,\infty)}(t).\]
For $f\in\mathcal S(\R)$, our Fourier transform convention gives
\[\widehat{Hf}(\xi)=\frac12\widehat f(\xi)-\frac i2\mathcal H\widehat f(\xi).\]
Since $\varpi_b\in A_2(\R)$, the weighted Hilbert transform estimate gives
\[\|\mathcal H\widehat f\|_{L^2(\varpi_b)}\lesssim_b\|\widehat f\|_{L^2(\varpi_b)}.\]
Consequently,
\[\|Hf\|_{H^b(\R)}=\left\|\frac12\widehat f-\frac i2\mathcal H\widehat f\right\|_{L^2(\varpi_b)}\lesssim_b\|\widehat f\|_{L^2(\varpi_b)}=\|f\|_{H^b(\R)}.\]
By density, this estimate extends to every $f\in H^b(\R)$.

For $a\in\R$, define $(\tau_af)(t):=f(t+a)$. Since $\tau_a$ is an isometry on $H^b(\R)$ and
\[\mathbf 1_{(a,\infty)}f=\tau_{-a}\left(\mathbf 1_{(0,\infty)}\tau_af\right),\]
we obtain
\[\|\mathbf 1_{(a,\infty)}f\|_{H^b(\R)}\lesssim_b\|f\|_{H^b(\R)}\]
with a constant independent of $a$. If $I=(a,c)$, then
\[\mathbf 1_I=\mathbf 1_{(a,\infty)}-\mathbf 1_{(c,\infty)},\]
and therefore
\[\|\mathbf 1_If\|_{H^b(\R)}\lesssim_b\|f\|_{H^b(\R)}\]
uniformly over all bounded intervals $I$.

We now apply this estimate mode by mode. For each $k\in\Z$, define
\[g_k(t):=e^{-it\omega(k)}\widehat u(t,k).\]
Then
\[\|u\|_{X^{\sigma,s,b}}^2=\sum_{k\in\Z}e^{2\sigma|k|}\bra{k}^{2s}\|g_k\|_{H^b(\R)}^2.\]
Since multiplication by $\mathbf 1_I$ commutes with the modulation factor $e^{-it\omega(k)}$, we also have
\[\|\mathbf 1_Iu\|_{X^{\sigma,s,b}}^2=\sum_{k\in\Z}e^{2\sigma|k|}\bra{k}^{2s}\|\mathbf 1_Ig_k\|_{H^b(\R)}^2.\]
Applying \eqref{eq:interval Hb multiplier} to each $g_k$ and summing over $k$, we obtain
\[\|\mathbf 1_Iu\|_{X^{\sigma,s,b}}^2\lesssim_b\sum_{k\in\Z}e^{2\sigma|k|}\bra{k}^{2s}\|g_k\|_{H^b(\R)}^2=\|u\|_{X^{\sigma,s,b}}^2.\]
This completes the proof.
\end{proof}

\begin{lem}\label{lem:strichartz embed}
Let $\alpha\ge2$ and
\[b>\frac{\alpha+2}{4(\alpha+1)}.\]
Then
\[\|u\|_{L_{t,x}^4(\R\times\T)}\lesssim_{\alpha,b}\|u\|_{X^{0,b}}.\]
\end{lem}

\begin{proof}
The proof follows the argument of Takaoka and Tsutsumi \cite{TakaokaTsutsumi2004}, which was inspired by Bourgain's proof of the periodic $L^4$ estimate \cite{Bourgain1993-1,Bourgain1993-2} and developed further by Tao \cite{Tao2001}. Related higher-order arguments can be found in \cite{Kwak2018,CKV2022}.

By the nesting property of the $X^{s,b}$ spaces, it suffices to consider
\begin{equation}\label{eq:b range}
\frac14<\frac{\alpha+2}{4(\alpha+1)}<b<\frac12.
\end{equation}

\medskip

Set
\[p:=\alpha+1,\qquad q:=4b-1.\]
Then $p\ge3$ and $q>1/p$.

Let $u=u_1+u_2$, where
\[\widehat u_1(t,k)=0\qquad\text{for }|k|>1.\]
Since
\[|u|^2\le2|u_1|^2+2|u_2|^2,\]
it suffices to estimate $\|u_1^2\|_{L^2_{t,x}}$ and $\|u_2^2\|_{L^2_{t,x}}$ separately.

We first consider $u_1$. By Plancherel's theorem,
\[\|u_1^2\|_{L_{t,x}^2}^2\lesssim\sum_{|k|\le2}\int_{\R}\left|\sum_{\substack{|k_1|\le1\\|k-k_1|\le1}}\int_{\R}|\widetilde u_1(\tau_1,k_1)||\widetilde u_1(\tau-\tau_1,k-k_1)|\; d\tau_1\right|^2 \; d\tau.\]
Define
\[\widetilde F_1(\tau,k):=\bra{\tau-\omega(k)}^b|\widetilde u_1(\tau,k)|.\]
We use the convolution estimate
\begin{equation}\label{eq:convolution estimate}
\int_{\R}\bra{a}^{-2b}\bra{c-a}^{-2b}\; da\lesssim_b\bra{c}^{1-4b},
\end{equation}
which holds for $1/4<b<1/2$. Since $b$ satisfies \eqref{eq:b range}, by \eqref{eq:convolution estimate},
\[\int_{\R}\frac{d\tau_1}{\bra{\tau_1-\omega(k_1)}^{2b}\bra{\tau-\tau_1-\omega(k-k_1)}^{2b}}\lesssim_b\bra{\tau-\omega(k_1)-\omega(k-k_1)}^{1-4b}\le1.\]
Therefore, the Cauchy--Schwarz and Minkowski inequalities give
\[\begin{aligned}
\|u_1^2\|_{L_{t,x}^2}^2
&\lesssim_b\sum_{|k|\le2}\left(\sum_{\substack{|k_1|\le1\\|k-k_1|\le1}}\left(\int_{\R^2}|\widetilde F_1(\tau_1,k_1)|^2|\widetilde F_1(\tau-\tau_1,k-k_1)|^2\; d\tau_1 \; d\tau\right)^{1/2}\right)^2\\
&\lesssim_b\|u_1\|_{X^{0,b}}^4.
\end{aligned}\]

We next consider $u_2$. Decompose
\[u_2=u_{2,+}+u_{2,-},\]
where $u_{2,+}$ and $u_{2,-}$ have spatial Fourier supports in $\{k>1\}$ and $\{k<-1\}$, respectively. Since
\[|u_2|^2\le2|u_{2,+}|^2+2|u_{2,-}|^2,\]
it suffices to estimate these two terms separately. We first consider $u_{2,+}$ and write it simply as $u_+$. Note that $\omega(k) = k^p$ for positive $k$.

By Plancherel's theorem,
\[\|u_+^2\|_{L_{t,x}^2}^2\lesssim\sum_{k>2}\int_{\R}\left|\sum_{\substack{k_1>1\\k-k_1>1}}\int_{\R}|\widetilde u_+(\tau_1,k_1)||\widetilde u_+(\tau-\tau_1,k-k_1)|\; d\tau_1\right|^2 \; d\tau.\]
Define
\[\widetilde F_+(\tau,k):=\bra{\tau-k^p}^b|\widetilde u_+(\tau,k)|.\]
By the Cauchy--Schwarz inequality,
\[\|u_+^2\|_{L_{t,x}^2}^2\lesssim M\|u_+\|_{X^{0,b}}^4,\]
where
\[M:=\sup_{\substack{\tau\in\R\\k>2}}\sum_{\substack{k_1>1\\k-k_1>1}}\int_{\R}\frac{d\tau_1}{\bra{\tau_1-k_1^p}^{2b}\bra{\tau-\tau_1-(k-k_1)^p}^{2b}}.\]
It follows from \eqref{eq:convolution estimate} that
\[M\lesssim_b\sup_{\substack{\tau\in\R\\k>2}}\sum_{\substack{k_1>1\\k-k_1>1}}\bra{\tau-k_1^p-(k-k_1)^p}^{-q}.\]

Fix $\tau\in\R$ and $k>2$, and define
\[h(x):=x^p+(k-x)^p-\tau,\qquad 0\le x\le k.\]
The function $h$ is symmetric about $x=k/2$, decreasing on $[0,k/2]$, and increasing on $[k/2,k]$. Moreover, since $p\ge3$, the function $r\mapsto r^{p-2}$ is convex on $[0,\infty)$. Hence, by Jensen's inequality,
\[x^{p-2}+(k-x)^{p-2}\ge2\left(\frac{x+(k-x)}2\right)^{p-2}=2^{3-p}k^{p-2}.\]
Therefore,
\[h''(x)=p(p-1)\left(x^{p-2}+(k-x)^{p-2}\right)\ge2^{3-p}p(p-1)k^{p-2}\gtrsim_pk^{p-2}.\]

For $R\ge1$, consider
\[E_R:=\{x\in[0,k]:|h(x)|\le R\}.\]
The intersection of $E_R$ with each of the intervals $[0,k/2]$ and $[k/2,k]$ is an interval. Let $[a,a+\ell]\subset[k/2,k]$ be one of these intervals. Since $|h(a)|\le R$, $|h(a+\ell)|\le R$, and $h$ is increasing on $[k/2,k]$, we have
\[0\le h(a+\ell)-h(a)\le2R.\]
By the fundamental theorem of calculus and integration by parts,
\[h(a+\ell)-h(a) = \int_0^\ell h'(a+r)\; dr = h'(a)\ell+\int_0^\ell(\ell-r)h''(a+r)\; dr.\]
Since $h'(a)\ge0$ and $h''(x)\gtrsim_pk^{p-2}$ on $[0,k]$, it follows that
\[2R\ge h(a+\ell)-h(a)\gtrsim_pk^{p-2}\int_0^\ell(\ell-r)\; dr\sim_pk^{p-2}\ell^2.\]
The same estimate holds on $[0,k/2]$ by symmetry. Hence,
\[\#\{k_1\in\Z:1<k_1<k-1,\ |h(k_1)|\le R\}\lesssim_p1+\min\left\{k,R^{1/2}k^{-(p-2)/2}\right\}.\]
Optimizing the minimum with respect to $k$, with the two terms balanced at $k=R^{1/p}$, yields
\[\min\left\{k,R^{1/2}k^{-(p-2)/2}\right\}\le R^{1/p}.\]
Therefore,
\[\#\{k_1\in\Z:1<k_1<k-1,\ |h(k_1)|\le R\}\lesssim_p1+R^{1/p}.\]

Using a dyadic decomposition, we obtain
\[\begin{aligned}
\sum_{\substack{k_1>1\\k-k_1>1}}\bra{h(k_1)}^{-q}
&\lesssim\sum_{m=0}^\infty2^{-mq}\#\{k_1:|h(k_1)|\le2^{m+1}\}\\
&\lesssim_p\sum_{m=0}^\infty2^{-mq}\left(1+2^{m/p}\right)\\
&\lesssim_{p,q}1,
\end{aligned}\]
where the last series converges since $q>1/p$. Thus,
\[M\lesssim_{\alpha,b}1,\]
and consequently
\[\|u_+^2\|_{L_{t,x}^2}^2\lesssim_{\alpha,b}\|u_+\|_{X^{0,b}}^4.\]

For the negative-frequency contribution, define
\[v(t,x):=u_{2,-}(-t,-x).\]
Then
\[\widetilde v(\tau,k)=\widetilde u_{2,-}(-\tau,-k),\]
so $v$ has spatial Fourier support in $\{k>1\}$. Since $\omega$ is odd,
\[\bra{\tau-\omega(k)}=\bra{-\tau-\omega(-k)}.\]
Therefore, the change of variables $(\tau,k)\mapsto(-\tau,-k)$ gives
\[\|v\|_{X^{0,b}}=\|u_{2,-}\|_{X^{0,b}}.\]
On the other hand, the invariance of the $L^4_{t,x}$ norm under reflections gives
\[\|v\|_{L_{t,x}^4}=\|u_{2,-}\|_{L_{t,x}^4}.\]
Applying the estimate for positive frequencies to $v$, we obtain
\[\|u_{2,-}^2\|_{L_{t,x}^2}^2\lesssim_{\alpha,b}\|u_{2,-}\|_{X^{0,b}}^4.\]
Combining the estimates for $u_1$, $u_{2,+}$, and $u_{2,-}$ yields
\[\|u\|_{L_{t,x}^4}\lesssim_{\alpha,b}\|u\|_{X^{0,b}}.\]
\end{proof}

\begin{lem}[Localized quartic estimate]\label{lem:localized quartic}
Let $\alpha\ge2$, $1/3<b<1/2$, and let $I\subset\R$ be a bounded interval. Suppose that $u_j\in X_I^{0,b}$ for $j=1,\ldots,4$, and let $M:\Z^4\to\mathbb C$ satisfy
\[\|M\|_{\ell^\infty}:=\sup_{k_1+k_2+k_3+k_4=0}|M(k_1,k_2,k_3,k_4)|<\infty.\]
Then
\[\left|\int_I\sum_{k_1+k_2+k_3+k_4=0}M(k_1,k_2,k_3,k_4)\prod_{j=1}^4\widehat u_j(t,k_j)\; dt\right|\lesssim_{\alpha,b}\|M\|_{\ell^\infty}\prod_{j=1}^4\|u_j\|_{X_I^{0,b}}.\]
\end{lem}

\begin{proof}
Fix $\varepsilon>0$. For each $j$, choose an extension $u_j^\ast\in X^{0,b}$ such that
\[\|u_j^\ast\|_{X^{0,b}}\le\|u_j\|_{X_I^{0,b}}+\varepsilon.\]
Define
\[v_j:=\mathbf 1_Iu_j^\ast.\]
Then $v_j=u_j$ on $I$ and $v_j=0$ outside $I$. By Lemma \ref{lem:interval truncation},
\[\|v_j\|_{X^{0,b}}\lesssim_b \|u_j^\ast\|_{X^{0,b}} \lesssim_b \|u_j\|_{X_I^{0,b}}+\varepsilon.\]

Set
\[\Lambda_I:=\int_I\sum_{k_1+k_2+k_3+k_4=0}M(k_1,k_2,k_3,k_4)\prod_{j=1}^4\widehat u_j(t,k_j)\; dt.\]
Since each $v_j$ vanishes outside $I$,
\[\Lambda_I=\int_{\R}\sum_{k_1+k_2+k_3+k_4=0}M(k_1,k_2,k_3,k_4)\prod_{j=1}^4\widehat v_j(t,k_j)\; dt.\]
Applying Fourier inversion in time, we obtain
\[\Lambda_I=\frac1{(2\pi)^3}\sum_{k_1+k_2+k_3+k_4=0}\int_{\tau_1+\tau_2+\tau_3+\tau_4=0}M(k_1,k_2,k_3,k_4)\prod_{j=1}^4\widetilde v_j(\tau_j,k_j)\; d\tau_1 \; d\tau_2 \; d\tau_3.\]
Therefore,
\[\begin{aligned}|\Lambda_I|&\lesssim\|M\|_{\ell^\infty}\sum_{k_1+k_2+k_3+k_4=0}\int_{\tau_1+\tau_2+\tau_3+\tau_4=0}\prod_{j=1}^4|\widetilde v_j(\tau_j,k_j)|\; d\tau_1 \; d\tau_2 \; d\tau_3.\end{aligned}\]

For each $j$, define $z_j$ by
\[\widetilde z_j(\tau,k):=|\widetilde v_j(\tau,k)|.\]
By space-time Fourier inversion, the convolution expression on the right-hand side equals, up to a harmless constant,
\[\int_{\R}\int_{\T}z_1(t,x)z_2(t,x)z_3(t,x)z_4(t,x)\; dx \; dt.\]
Hence, H\"older's inequality and Lemma \ref{lem:strichartz embed} give
\[|\Lambda_I| \lesssim_{\alpha,b} \|M\|_{\ell^\infty}\prod_{j=1}^4\|v_j\|_{X^{0,b}} \lesssim_b \|M\|_{\ell^\infty}\prod_{j=1}^4\left(\|u_j\|_{X_I^{0,b}}+\varepsilon\right).\]
Letting $\varepsilon \to 0$ completes the proof.
\end{proof}

\section{Normal form reduction}\label{sec:normal form}
Throughout the remainder of the paper, we work under the mean-zero assumption. The general case is recovered by the Galilean transformation in Remark \ref{rem:recover mean}. 

In this section, we consider smooth real-valued solutions to \eqref{eq:dgkdv}. We first introduce the interaction representation
\[v(t):=e^{-it\omega(D_x)}u(t),\]
where $\omega(D_x)$ denotes the Fourier multiplier with symbol $\omega(k)$. Equivalently,
\begin{equation}\label{eq:interaction representation}
\hat v(t,k)=e^{-it\omega(k)}\hat u(t,k).
\end{equation}
Taking the spatial Fourier transform of \eqref{eq:dgkdv}, we obtain
\begin{equation}\label{eq:u fourier evolution}
\partial_t\hat u(t,k)-i\omega(k)\hat u(t,k)=-\frac{ik}{2}\sum_{k_1+k_2=k}\hat u(t,k_1)\hat u(t,k_2).
\end{equation}
For $k_1,k_2\in\mathbb Z$, define the quadratic interaction phase by
\[\Omega(k_1,k_2):=\omega(k_1)+\omega(k_2)-\omega(k_1+k_2).\]
Using \eqref{eq:interaction representation} in \eqref{eq:u fourier evolution}, we obtain
\begin{equation}\label{eq:v evolution}
\partial_t\hat v(t,k)=-\frac{ik}{2}\sum_{k_1+k_2=k}e^{it\Omega(k_1,k_2)}\hat v(t,k_1)\hat v(t,k_2).
\end{equation}

We next define the Gevrey profile
\[V(t):=e^{\sigma|D_x|}v(t)=e^{\sigma|D_x|}e^{-it\omega(D_x)}u(t).\]
In terms of the Fourier coefficients,
\begin{equation}\label{eq:Gevrey profile}
\hat V(t,k)=e^{\sigma|k|}\hat v(t,k)=e^{\sigma|k|}e^{-it\omega(k)}\hat u(t,k).
\end{equation}
Multiplying \eqref{eq:v evolution} by $e^{\sigma|k|}$ gives
\begin{equation}\label{eq:V evolution}
\partial_t\hat V(t,k)=-\frac{ik}{2}\sum_{k_1+k_2=k}e^{\sigma(|k|-|k_1|-|k_2|)}e^{it\Omega(k_1,k_2)}\hat V(t,k_1)\hat V(t,k_2).
\end{equation}
Since $u$ is real-valued and has zero spatial mean, we have
\begin{equation}\label{eq:V reality mean}
\hat V(t,-k)=\overline{\hat V(t,k)},\qquad \hat V(t,0)=0.
\end{equation}

For the cubic interaction, define
\[\Gamma_3:=\{(k_1,k_2,k_3)\in(\mathbb Z\setminus\{0\})^3:k_1+k_2+k_3=0\}.\]
Define the cubic phase function by
\[\Phi_3(k_1,k_2,k_3):=\sum_{j=1}^3\omega(k_j).\]
If $k_1+k_2+k_3=0$, then
\begin{equation}\label{eq:quadratic cubic phase relation}
\Omega(k_1,k_2)=\Phi_3(k_1,k_2,k_3).
\end{equation}

\begin{lem}[Factorization of the cubic phase]\label{lem:phase factorization}
Let $\alpha\ge2$ and $(k_1,k_2,k_3)\in\Gamma_3$. Then
\begin{equation}\label{eq:phase factorization}
\Phi_3(k_1,k_2,k_3)=k_1k_2k_3Q_3(k_1,k_2,k_3),
\end{equation}
where $Q_3$ is a real-valued symmetric function satisfying
\begin{equation}\label{eq:Q3 bound}
|Q_3(k_1,k_2,k_3)|\sim_\alpha K_{\max}^{\alpha-2},
\end{equation}
with
\[K_{\max}:=\max\{|k_1|,|k_2|,|k_3|\}.\]
In particular,
\[\Phi_3(k_1,k_2,k_3)\ne0\]
for every $(k_1,k_2,k_3)\in\Gamma_3$.
\end{lem}

\begin{proof}
When $\alpha=2j$, $j\in\mathbb N$, the factorization \eqref{eq:phase factorization} and the bound \eqref{eq:Q3 bound} coincide with Lemma 2.1(a) in \cite{HK}. We provide a direct proof for every real $\alpha\ge2$.

Define
\[Q_3(k_1,k_2,k_3):=\frac{\Phi_3(k_1,k_2,k_3)}{k_1k_2k_3}.\]
This is well-defined on $\Gamma_3$ since $k_1k_2k_3\ne0$. Both $\Phi_3(k_1,k_2,k_3)$ and $k_1k_2k_3$ are symmetric under permutations of $(k_1,k_2,k_3)$. Hence, $Q_3$ is also symmetric under permutations of the frequencies.

Since $(k_1,k_2,k_3)\in\Gamma_3$, we have
\[k_1+k_2+k_3=0\]
and $k_j\ne0$ for each $j$. Note that the three frequencies cannot have the same sign. Therefore, exactly two of them have the same sign and the remaining one has the opposite sign. In particular, even if the original $k_1$ and $k_2$ have different signs, we may permute the indices so that the two frequencies having the same sign are relabeled as $k_1$ and $k_2$.

After this permutation, either $k_1,k_2>0$ or $k_1,k_2<0$. Since $\omega$ is odd, we have
\[\Phi_3(-k_1,-k_2,-k_3)=-\Phi_3(k_1,k_2,k_3).\]
Moreover,
\[(-k_1)(-k_2)(-k_3)=-k_1k_2k_3.\]
Consequently,
\[Q_3(-k_1,-k_2,-k_3)=Q_3(k_1,k_2,k_3).\]
Thus, it is enough to consider
\[k_1=a,\qquad k_2=b,\qquad k_3=-(a+b),\]
where $a,b>0$.

In this configuration,
\[\Phi_3(k_1,k_2,k_3)=a^{\alpha+1}+b^{\alpha+1}-(a+b)^{\alpha+1}\]
and
\[k_1k_2k_3=-ab(a+b).\]
Hence,
\[Q_3(k_1,k_2,k_3)=\frac{(a+b)^{\alpha+1}-a^{\alpha+1}-b^{\alpha+1}}{ab(a+b)}.\]

Set $p:=\alpha+1$. By the fundamental theorem of calculus,
\[(a+b)^p-a^p-b^p=p(p-1)\int_0^b\int_0^a(x+y)^{p-2}\; dx \; dy.\]
Since $x+y\le a+b$ on $[0,a]\times[0,b]$, we obtain
\[(a+b)^p-a^p-b^p\lesssim_p ab(a+b)^{p-2}.\]
For the reverse estimate, we restrict the integration to $[a/2,a]\times[b/2,b]$. On this rectangle,
\[x+y\ge\frac{a+b}{2},\]
and therefore
\[(a+b)^p-a^p-b^p\gtrsim_p ab(a+b)^{p-2}.\]
It follows that
\[(a+b)^p-a^p-b^p\sim_p ab(a+b)^{p-2}.\]
Consequently,
\[|Q_3(k_1,k_2,k_3)|\sim_\alpha(a+b)^{p-3}=(a+b)^{\alpha-2}.\]
Since $a,b>0$ and $k_3=-(a+b)$, we have
\[K_{\max}=a+b.\]
Thus,
\[|Q_3(k_1,k_2,k_3)|\sim_\alpha K_{\max}^{\alpha-2}.\]
This completes the proof.
\end{proof}

Define the Gevrey energy by
\begin{equation}\label{eq:Gevrey energy}
E_\sigma(t):=\sum_{k\in\mathbb Z}|\hat V(t,k)|^2.
\end{equation}
By \eqref{eq:Gevrey profile},
\[E_\sigma(t)=\sum_{k\in\mathbb Z}e^{2\sigma|k|}|\hat u(t,k)|^2=\|u(t)\|_{G^{\sigma,0}}^2.\]
Differentiating \eqref{eq:Gevrey energy}, we obtain
\[\frac d{dt}E_\sigma(t)=2\Re\sum_{k\in\mathbb Z}\overline{\hat V(t,k)}\partial_t\hat V(t,k).\]
Substituting \eqref{eq:V evolution} gives
\[\frac d{dt}E_\sigma(t)=-\Re\sum_{k\in\mathbb Z}\sum_{k_1+k_2=k} ike^{\sigma(|k|-|k_1|-|k_2|)}e^{it\Omega(k_1,k_2)} \overline{\hat V(t,k)}\hat V(t,k_1)\hat V(t,k_2).\]
Using $\overline{\hat V(t,k)}=\hat V(t,-k)$ and setting $k_3:=-k$, we have $k_1+k_2+k_3=0$. By \eqref{eq:quadratic cubic phase relation},
\[\Omega(k_1,k_2)=\Phi_3(k_1,k_2,k_3).\]
All terms involving a zero frequency vanish by \eqref{eq:V reality mean}. Therefore,
\begin{equation}\label{eq:energy unsymmetrised}
\begin{aligned}
\frac d{dt}E_\sigma(t)=\Re\sum_{(k_1,k_2,k_3)\in\Gamma_3} &ik_3e^{\sigma(|k_3|-|k_1|-|k_2|)}e^{it\Phi_3(k_1,k_2,k_3)}\\
&\times\hat V(t,k_1)\hat V(t,k_2)\hat V(t,k_3).
\end{aligned}
\end{equation}

The summation set $\Gamma_3$, the phase $\Phi_3$, and the product of the Fourier coefficients are invariant under permutations of $(k_1,k_2,k_3)$. Symmetrizing the multiplier in \eqref{eq:energy unsymmetrised}, we define
\begin{equation}\label{eq:M3 def}
M_3(k_1,k_2,k_3):=\frac i3\sum_{j=1}^3k_j\exp\left(\sigma\left(|k_j|-\sum_{\substack{1\le m\le3\\m\ne j}}|k_m|\right)\right).
\end{equation}
Moreover, since $k_1+k_2+k_3=0$, the multiplier \eqref{eq:M3 def} can equivalently be written as
\begin{equation}\label{eq:M3 commutator}
M_3(k_1,k_2,k_3)=\frac i3\sum_{j=1}^3k_j\left[\exp\left(\sigma\left(|k_j|-\sum_{\substack{1\le m\le3\\m\ne j}}|k_m|\right)\right)-1\right].
\end{equation}
Then
\begin{equation}\label{eq:energy profile}
\frac{d}{dt}E_\sigma(t)=\Re\sum_{(k_1,k_2,k_3)\in\Gamma_3}M_3(k_1,k_2,k_3)e^{it\Phi_3(k_1,k_2,k_3)}\prod_{j=1}^3\hat V(t,k_j).
\end{equation}
We suppress the dependence of $M_3$ on $\sigma$ in the notation.	

By Lemma \ref{lem:phase factorization}, $\Phi_3$ does not vanish on $\Gamma_3$. Hence,
\begin{equation}\label{eq:oscillation derivative}
e^{it\Phi_3}=-\frac i{\Phi_3}\partial_t(e^{it\Phi_3}).
\end{equation}
Define the cubic boundary term by
\begin{equation}\label{eq:B3 def}
B(t):=\operatorname{Im}\sum_{(k_1,k_2,k_3)\in\Gamma_3}\frac{M_3(k_1,k_2,k_3)}{\Phi_3(k_1,k_2,k_3)}e^{it\Phi_3(k_1,k_2,k_3)}\prod_{j=1}^3\hat V(t,k_j).
\end{equation}
Substituting \eqref{eq:oscillation derivative} into \eqref{eq:energy profile} and applying the product rule, we obtain
\begin{equation}\label{eq:integration by parts}
\frac{d}{dt}E_\sigma(t) =\frac{d}{dt}B(t) -\Im\sum_{(k_1,k_2,k_3)\in\Gamma_3} \frac{M_3(k_1,k_2,k_3)}{\Phi_3(k_1,k_2,k_3)} e^{it\Phi_3(k_1,k_2,k_3)}\partial_t\left(\prod_{j=1}^3\hat V(t,k_j)\right).
\end{equation}

The multiplier
\[\frac{M_3(k_1,k_2,k_3)}{\Phi_3(k_1,k_2,k_3)}e^{it\Phi_3(k_1,k_2,k_3)}\]
is symmetric in $(k_1,k_2,k_3)$. Therefore, the three terms obtained by differentiating the cubic product give the same contribution after relabeling the summation indices. It follows that
\begin{equation}\label{eq:cubic product derivative}
\begin{aligned}
&\sum_{(k_1,k_2,k_3)\in\Gamma_3}\frac{M_3(k_1,k_2,k_3)}{\Phi_3(k_1,k_2,k_3)}e^{it\Phi_3(k_1,k_2,k_3)}\partial_t\left(\prod_{j=1}^3\hat V(t,k_j)\right)\\
=&~{}3\sum_{(\ell,k_3,k_4)\in\Gamma_3}\frac{M_3(\ell,k_3,k_4)}{\Phi_3(\ell,k_3,k_4)}e^{it\Phi_3(\ell,k_3,k_4)}\partial_t\hat V(t,\ell)\hat V(t,k_3)\hat V(t,k_4).
\end{aligned}
\end{equation}

For the quartic interaction, define
\[\Gamma_4:=\{(k_1,k_2,k_3,k_4)\in(\mathbb Z\setminus\{0\})^4:k_1+k_2+k_3+k_4=0\}\]
and
\[\Gamma_4^\ast:=\{(k_1,k_2,k_3,k_4)\in\Gamma_4:k_1+k_2\ne0\}.\]
Define the quartic phase function by
\[\Phi_4(k_1,k_2,k_3,k_4):=\sum_{j=1}^4\omega(k_j).\]
By \eqref{eq:V evolution},
\begin{equation}\label{eq:new profile}
\partial_t\hat V(t,\ell)=-\frac{i\ell}{2}\sum_{k_1+k_2=\ell}e^{\sigma(|\ell|-|k_1|-|k_2|)}e^{it\Omega(k_1,k_2)}\hat V(t,k_1)\hat V(t,k_2).
\end{equation}
For $\ell=k_1+k_2$, the phases satisfy
\[\begin{aligned}\Phi_3(\ell,k_3,k_4)+\Omega(k_1,k_2)&=\omega(\ell)+\omega(k_3)+\omega(k_4)+\omega(k_1)+\omega(k_2)-\omega(\ell)\\&=\Phi_4(k_1,k_2,k_3,k_4).\end{aligned}\]
Since $(\ell,k_3,k_4)\in\Gamma_3$, we have $\ell\ne0$. Terms for which $k_1=0$ or $k_2=0$ vanish by \eqref{eq:V reality mean}. Hence, after setting
\[k_{12}:=k_1+k_2,\]
the resulting quartic summation is over $\Gamma_4^\ast$.

Substituting \eqref{eq:new profile} into \eqref{eq:cubic product derivative}, the second term on the right-hand side of \eqref{eq:integration by parts} becomes
\[-\Im\left(-\frac{3i}{2}\sum_{(k_1,k_2,k_3,k_4)\in\Gamma_4^\ast}\frac{M_3(k_{12},k_3,k_4)}{\Phi_3(k_{12},k_3,k_4)}k_{12}e^{\sigma(|k_{12}|-|k_1|-|k_2|)}e^{it\Phi_4(k_1,k_2,k_3,k_4)}\prod_{j=1}^4\hat V(t,k_j)\right).\]
Converting the imaginary part into a real part through the factor $-i$, we define the quartic remainder by
\begin{equation}\label{eq:R4 def}
\begin{aligned}
R_4(t):=\frac32\Re\sum_{(k_1,k_2,k_3,k_4)\in\Gamma_4^\ast}\frac{M_3(k_{12},k_3,k_4)}{\Phi_3(k_{12},k_3,k_4)}k_{12}e^{\sigma(|k_{12}|-|k_1|-|k_2|)}e^{it\Phi_4(k_1,k_2,k_3,k_4)}\prod_{j=1}^4\hat V(t,k_j).
\end{aligned}
\end{equation}
We have therefore obtained
\begin{equation}\label{eq:modified energy evolution}
\frac{d}{dt}E_\sigma(t)=\frac d{dt}B(t)+R_4(t).
\end{equation}
Integrating \eqref{eq:modified energy evolution} over $[0,t]$, we obtain the following normal form identity.

\begin{prop}\label{prop:normal form identity}
Let $V$ be a smooth solution to \eqref{eq:V evolution} satisfying \eqref{eq:V reality mean}. Then, for every $t$ in the lifespan of $V$,
\[E_\sigma(t) = E_\sigma(0) + B(t) - B(0)+\int_0^tR_4(t')\; dt',\]
where $E_\sigma$, $B$, and $R_4$ are defined by \eqref{eq:Gevrey energy}, \eqref{eq:B3 def}, and \eqref{eq:R4 def}, respectively.
\end{prop}

\section{Local-in-time almost conservation laws}\label{sec:estimates}
Throughout this section, $u$ denotes a smooth real-valued solution to \eqref{eq:dgkdv}, and $E_\sigma$, $B$, and $R_4$ are defined as in Section \ref{sec:normal form}.

\begin{lem}
Let $\alpha\ge2$ and $\sigma\ge0$. For every $(k_1,k_2,k_3)\in\Gamma_3$, the multiplier associated with the boundary term satisfies
\begin{equation}\label{eq:symbol B}
\left|\frac{M_3(k_1,k_2,k_3)}{\Phi_3(k_1,k_2,k_3)}\right|\lesssim_\alpha\sigma K_{\max}^{1-\alpha},
\end{equation}
where
\[K_{\max}:=\max\{|k_1|,|k_2|,|k_3|\}.\]
Furthermore, for every $(k_1,k_2,k_3,k_4)\in\Gamma_4^\ast$, writing $k_{12}:=k_1+k_2$, we have
\begin{equation}\label{eq:symbol R}
\left|\frac{M_3(k_{12},k_3,k_4)}{\Phi_3(k_{12},k_3,k_4)}k_{12}\right|\lesssim_\alpha\sigma.
\end{equation}
Consequently,
\begin{equation}\label{eq:full symbol R}
\left|\frac{M_3(k_{12},k_3,k_4)}{\Phi_3(k_{12},k_3,k_4)}k_{12}e^{\sigma(|k_{12}|-|k_1|-|k_2|)}\right|\lesssim_\alpha\sigma.
\end{equation}
\end{lem}

\begin{proof}
We first prove \eqref{eq:symbol B}. Let
\[K_{\max}>K_{\mathrm{med}}\ge K_{\min}>0\]
be the numbers $|k_1|$, $|k_2|$, and $|k_3|$ arranged in decreasing order. Since $k_1+k_2+k_3=0$ and all three frequencies are nonzero, the two frequencies having magnitudes $K_{\mathrm{med}}$ and $K_{\min}$ have the same sign, while the remaining frequency has the opposite sign. It follows that
\[K_{\max}=K_{\mathrm{med}}+K_{\min}.\]

Recall the commutator representation \eqref{eq:M3 commutator},
\[M_3(k_1,k_2,k_3)=\frac i3\sum_{j=1}^3k_j\left[\exp\left(\sigma\left(|k_j|-\sum_{\substack{1\le m\le3\\m\ne j}}|k_m|\right)\right)-1\right].\]
For the largest frequency, the exponent vanishes since
\[K_{\max}-K_{\mathrm{med}}-K_{\min}=0.\]
The corresponding contribution to $M_3$ is therefore zero. For the remaining two frequencies, the exponents are
\[-2K_{\min}\qquad\text{and}\qquad-2K_{\mathrm{med}},\]
respectively. Using $1-e^{-x}\le x$ for $x\ge0$, we obtain
\[\begin{aligned}
|M_3(k_1,k_2,k_3)|
&\le\frac13K_{\mathrm{med}}\left(1-e^{-2\sigma K_{\min}}\right)+\frac13K_{\min}\left(1-e^{-2\sigma K_{\mathrm{med}}}\right)\\
&\le\frac{2\sigma}{3}K_{\mathrm{med}}K_{\min}+\frac{2\sigma}{3}K_{\min}K_{\mathrm{med}}\\
&\lesssim\sigma K_{\mathrm{med}}K_{\min}.
\end{aligned}\]
By Lemma \ref{lem:phase factorization},
\[|\Phi_3(k_1,k_2,k_3)|\sim_\alpha|k_1k_2k_3|K_{\max}^{\alpha-2}=K_{\max}^{\alpha-1}K_{\mathrm{med}}K_{\min}.\]
Consequently,
\[\left|\frac{M_3(k_1,k_2,k_3)}{\Phi_3(k_1,k_2,k_3)}\right|\lesssim_\alpha\frac{\sigma K_{\mathrm{med}}K_{\min}}{K_{\max}^{\alpha-1}K_{\mathrm{med}}K_{\min}}=\sigma K_{\max}^{1-\alpha},\]
which proves \eqref{eq:symbol B}.

We next prove \eqref{eq:symbol R}. Let
\[K:=\max\{|k_{12}|,|k_3|,|k_4|\}.\]
Since $(k_1,k_2,k_3,k_4)\in\Gamma_4^\ast$, we have $k_{12}\ne0$. Moreover,
\[k_{12}+k_3+k_4=0,\]
and hence $(k_{12},k_3,k_4)\in\Gamma_3$. Applying \eqref{eq:symbol B} to this triple gives
\[\left|\frac{M_3(k_{12},k_3,k_4)}{\Phi_3(k_{12},k_3,k_4)}\right|\lesssim_\alpha\sigma K^{1-\alpha}.\]
It follows that
\[\left|\frac{M_3(k_{12},k_3,k_4)}{\Phi_3(k_{12},k_3,k_4)}k_{12}\right|\lesssim_\alpha\sigma|k_{12}|K^{1-\alpha}\le\sigma K^{2-\alpha}.\]
Since $\alpha\ge2$ and $K\ge1$, we have $K^{2-\alpha}\le1$. This proves \eqref{eq:symbol R}.

Finally, the triangle inequality gives
\[|k_{12}|\le|k_1|+|k_2|,\]
and therefore
\[e^{\sigma(|k_{12}|-|k_1|-|k_2|)}\le1.\]
Combining this with \eqref{eq:symbol R} proves \eqref{eq:full symbol R}.
\end{proof}

\begin{lem}\label{lem:boundary bound}
Let $\alpha\ge2$ and $\sigma\ge0$. Then the boundary term $B(t)$ defined in \eqref{eq:B3 def} satisfies
\begin{equation}\label{eq:boundary bound}
|B(t)|\lesssim_\alpha\sigma E_\sigma(t)^{3/2},
\end{equation}
where $E_\sigma(t)$ is defined as in \eqref{eq:Gevrey energy}.
\end{lem}

\begin{proof}
By the definition of $B(t)$ and \eqref{eq:symbol B},
\[\begin{aligned}
|B(t)|
&\le\sum_{(k_1,k_2,k_3)\in\Gamma_3}\left|\frac{M_3(k_1,k_2,k_3)}{\Phi_3(k_1,k_2,k_3)}\right|\prod_{j=1}^3|\hat V(t,k_j)|\\
&\lesssim_\alpha\sigma\sum_{(k_1,k_2,k_3)\in\Gamma_3}K_{\max}^{1-\alpha}\prod_{j=1}^3|\hat V(t,k_j)|.
\end{aligned}\]
Since $\alpha\ge2$ and $K_{\max}\ge1$,
\[K_{\max}^{1-\alpha}\le K_{\max}^{-1}\le|k_1|^{-1/3}|k_2|^{-1/3}|k_3|^{-1/3}.\]
Thus,
\[|B(t)|\lesssim_\alpha\sigma\sum_{(k_1,k_2,k_3)\in\Gamma_3}\prod_{j=1}^3\left(|k_j|^{-1/3}|\hat V(t,k_j)|\right).\]

Define $h(t)$ by
\[\hat h(t,0):=0,\qquad \hat h(t,k):=|k|^{-1/3}|\hat V(t,k)|\quad\text{for }k\ne0.\]
By Fourier inversion and the condition $k_1+k_2+k_3=0$,
\[\sum_{(k_1,k_2,k_3)\in\Gamma_3}\prod_{j=1}^3\hat h(t,k_j)\lesssim\int_{\mathbb T}|h(t,x)|^3\; dx.\]
Hence, by the Sobolev embedding $H^{1/6}(\mathbb T)\hookrightarrow L^3(\mathbb T)$,
\[|B(t)|\lesssim_\alpha\sigma\|h(t)\|_{L^3}^3\lesssim_\alpha\sigma\|h(t)\|_{H^{1/6}}^3.\]
Since $k\ne0$ on the support of $\hat h$,
\[\begin{aligned}
\|h(t)\|_{H^{1/6}}^2
&=\sum_{k\ne0}\langle k\rangle^{1/3}|k|^{-2/3}|\hat V(t,k)|^2\\
&\lesssim\sum_{k\ne0}|\hat V(t,k)|^2\\
&\le E_\sigma(t).
\end{aligned}\]
Therefore,
\[|B(t)|\lesssim_\alpha\sigma E_\sigma(t)^{3/2},\]
which proves \eqref{eq:boundary bound}.
\end{proof}

\begin{lem}\label{lem:remainder bound}
Let $\alpha\ge2$, $\sigma\ge0$, and $0<\delta\le1$. Then, for every $t\in[0,\delta]$, the time-integrated remainder satisfies
\begin{equation}\label{eq:remainder bound}
\left|\int_0^tR_4(t')\; dt'\right|\lesssim_\alpha\sigma\delta^{1/4}\|u\|_{Y_\delta^{\sigma,0}}^4.
\end{equation}
\end{lem}

\begin{proof}
Set
\[U(t):=e^{\sigma|D_x|}u(t).\]
By the definition of the Gevrey profile in \eqref{eq:Gevrey profile},
\[\hat V(t,k)=e^{-it\omega(k)}\hat U(t,k).\]
Therefore, on $\Gamma_4$,
\[e^{it\Phi_4(k_1,k_2,k_3,k_4)}\prod_{j=1}^4\hat V(t,k_j)=\prod_{j=1}^4\hat U(t,k_j).\]

For $(k_1,k_2,k_3,k_4)\in\Gamma_4^\ast$, define
\[M_4(k_1,k_2,k_3,k_4):=\frac{M_3(k_{12},k_3,k_4)}{\Phi_3(k_{12},k_3,k_4)}k_{12}e^{\sigma(|k_{12}|-|k_1|-|k_2|)},\]
where $k_{12}:=k_1+k_2$. By \eqref{eq:full symbol R},
\begin{equation}\label{eq:M4 bound}
|M_4(k_1,k_2,k_3,k_4)|\lesssim_\alpha\sigma.
\end{equation}
We extend $M_4$ by zero to all integer quadruples satisfying
\[k_1+k_2+k_3+k_4=0\]
but not belonging to $\Gamma_4^\ast$. It follows from \eqref{eq:R4 def} that
\[\int_0^tR_4(t')\; dt'=\frac32\Re\int_0^t\sum_{k_1+k_2+k_3+k_4=0}M_4(k_1,k_2,k_3,k_4)\prod_{j=1}^4\hat U(t',k_j)\; dt'.\]

The case $t=0$ is immediate. Let $0<t\le\delta$, and set
\[b_0:=\frac38,\qquad b_1:=\frac7{16}.\]
Then
\[\frac13<b_0<b_1<\frac12,\qquad4(b_1-b_0)=\frac14.\]
Applying Lemma \ref{lem:localized quartic} on the interval $[0,t]$ and using \eqref{eq:M4 bound}, we obtain
\begin{equation}\label{eq:localized remainder X}
\left|\int_0^tR_4(t')\; dt'\right|\lesssim_\alpha\sigma\|U\|_{X_{[0,t]}^{0,b_0}}^4.
\end{equation}

Choose an extension $u^\ast$ of $u$ from $[0,\delta]$ to $\mathbb R$ such that
\[\|u^\ast\|_{Y^{\sigma,0}}\le2\|u\|_{Y_\delta^{\sigma,0}},\]
and define
\[U^\ast:=e^{\sigma|D_x|}u^\ast.\]
Then
\[\|U^\ast\|_{Y^0}=\|u^\ast\|_{Y^{\sigma,0}}.\]

Fix $\eta\in C_c^\infty(\mathbb R)$ such that
\[\eta=1\text{ on }[0,1].\]
Define
\[F_t(t'):=\eta(t'/t)U^\ast(t').\]
Since $F_t=U$ on $[0,t]$, the definition of the restriction norm gives
\[\|U\|_{X_{[0,t]}^{0,b_0}}\le\|F_t\|_{X^{0,b_0}}.\]
By Lemma \ref{lem:time localization},
\[\begin{aligned}
\|F_t\|_{X^{0,b_0}}
&\lesssim t^{b_1-b_0}\|U^\ast\|_{X^{0,b_1}}\\
&\le\delta^{b_1-b_0}\|U^\ast\|_{Y^0}\\
&\lesssim\delta^{b_1-b_0}\|u\|_{Y_\delta^{\sigma,0}}.
\end{aligned}\]
Consequently,
\begin{equation}\label{eq:localized U bound}
\|U\|_{X_{[0,t]}^{0,b_0}}\lesssim\delta^{b_1-b_0}\|u\|_{Y_\delta^{\sigma,0}}.
\end{equation}
Combining \eqref{eq:localized remainder X} and \eqref{eq:localized U bound}, we conclude
\[\begin{aligned}
\left|\int_0^tR_4(t')\; dt'\right|
&\lesssim_\alpha\sigma\delta^{4(b_1-b_0)}\|u\|_{Y_\delta^{\sigma,0}}^4\\
&=\sigma\delta^{1/4}\|u\|_{Y_\delta^{\sigma,0}}^4.
\end{aligned}\]
This proves \eqref{eq:remainder bound}.
\end{proof}

Combining Proposition \ref{prop:normal form identity} with Lemmas \ref{lem:boundary bound} and \ref{lem:remainder bound} yields the corresponding local almost conservation estimate. In the proof of Theorem \ref{thm:main}, however, we keep the boundary and remainder terms separate, since only the latter accumulates over successive local intervals.


\section{Proof of Theorem \ref{thm:main}}\label{sec:proof}
In this section, we prove the global existence and quantitative persistence assertions in Theorem~\ref{thm:main}. The unconditional uniqueness and continuous-dependence assertions were justified in Remarks \ref{rem:unconditional uniqueness} and \ref{rem:continuous dependence}, respectively. 

The normal form identity and the estimates established in Sections \ref{sec:normal form} and \ref{sec:estimates} were derived for smooth solutions. By a standard approximation argument using sharp Fourier projections, they remain valid for the local solutions constructed in Theorem \ref{thm:hks lwp}. We therefore apply this identity and estimates below without further comment.

The equation is invariant under the transformation
\[u(t,x)\mapsto u(-t,-x).\]
It is therefore sufficient to establish the result for positive times.

\subsection{The case $s=0$}

Let $u_0\in G^{\sigma_0,0}(\mathbb T)$ and set
\[E_\ast:=E_{\sigma_0}(0)=\|u_0\|_{G^{\sigma_0,0}}^2.\]
The case $E_\ast=0$ is immediate, so we assume that $E_\ast>0$.

By Theorem \ref{thm:hks lwp}, there exist constants $a,c_0,C_0>0$, depending only on $\alpha$ and uniform with respect to the analytic radius, such that the following holds. If $\rho>0$ and
\[\|u(t_0)\|_{G^{\rho,0}}^2\le2E_\ast,\]
then the solution exists on $[t_0,t_0+\delta]$, where
\[\delta:=\min\left\{1,c_0\left(1+\sqrt{2E_\ast}\right)^{-a}\right\},\]
and satisfies
\begin{equation}\label{eq:uniform local bound}
\|u\|_{Y_{[t_0,t_0+\delta]}^{\rho,0}}^2\le C_{\mathrm{loc}}E_\ast,
\end{equation}
where
\[C_{\mathrm{loc}}:=2C_0^2.\]
Thus, $\delta$ depends only on $\alpha$ and $E_\ast$, while $C_{\mathrm{loc}}$ depends only on $\alpha$.

Let $C_B>0$ be the constant in Lemma \ref{lem:boundary bound}. By time translation, Lemma \ref{lem:remainder bound} and \eqref{eq:uniform local bound} imply that there exists $C_R>0$, depending only on $\alpha$, such that for every interval $J$ of length at most $\delta$,
\begin{equation}\label{eq:local remainder step}
\left|\int_JR_4(t')\; dt'\right|\le C_R\sigma\delta^{1/4}E_\ast^2
\end{equation}
whenever
\[\sup_{t\in J}E_\sigma(t)\le2E_\ast.\]

Choose
\begin{equation}\label{eq:analyticity constant}
c_\ast:=\min\left\{\frac{\delta}{4C_B(1+2^{3/2})E_\ast^{1/2}},\frac{\delta^{3/4}}{8C_RE_\ast}\right\}
\end{equation}
and, for a fixed $T>0$, set
\begin{equation}\label{eq:sigma choice}
\sigma:=\min\left\{\sigma_0,\frac{c_\ast}{T}\right\}.
\end{equation}

If $0<T\le\delta$, Theorem \ref{thm:hks lwp} gives
\[u\in C([-T,T];G^{\sigma_0,0}(\mathbb T))\subset C([-T,T];G^{\sigma,0}(\mathbb T)).\]
We thus assume that $T>\delta$. In particular,
\begin{equation}\label{eq:sigma boundary smallness}
\sigma\le\frac{c_\ast}{T}\le\frac{c_\ast}{\delta}\le\frac{1}{4C_B(1+2^{3/2})E_\ast^{1/2}}.
\end{equation}

Let $u$ denote the maximal forward solution by Theorem \ref{thm:hks lwp}. By the embedding
\[Y_I^{\sigma,0}\hookrightarrow C(I;G^{\sigma,0}(\mathbb T)),\]
the solution is continuous in $G^{\sigma,0}(\mathbb T)$ on every compact subinterval of its maximal lifespan. Consequently,
\[E_\sigma(t)=\|u(t)\|_{G^{\sigma,0}}^2\]
is continuous throughout the maximal lifespan.

Define
\[T^\ast:=\sup\left\{T_0\in[0,T]:u\text{ exists on }[0,T_0]\text{ and }\sup_{0\le t\le T_0}E_\sigma(t)\le2E_\ast\right\}.\]
Since
\[E_\sigma(0)\le E_{\sigma_0}(0)=E_\ast,\]
the continuity of $E_\sigma$ implies that $T^\ast>0$.

Fix $0<t<T^\ast$. By the definition of $T^\ast$, the solution exists on $[0,t]$ and
\[\sup_{0\le t'\le t}E_\sigma(t')\le2E_\ast.\]
We split $[0,t]$ into $N$ subintervals of length at most $\delta$. Since $T>\delta$, we may choose $N$ such that
\[N\le1+\frac{t}{\delta}\le\frac{2T}{\delta}.\]
For every interval $J$ in the partition, we have
\begin{equation}\label{eq:bootstrap assumption}
\sup_{t'\in J}E_\sigma(t')\le2E_\ast.
\end{equation}
Therefore, \eqref{eq:local remainder step} gives
\[\left|\int_JR_4(t')\; dt'\right|\le C_R\sigma\delta^{1/4}E_\ast^2.\]
Summing over the partition, we obtain
\[\left|\int_0^tR_4(t')\; dt'\right|\le NC_R\sigma\delta^{1/4}E_\ast^2\le2C_RT\sigma\delta^{-3/4}E_\ast^2.\]
Using \eqref{eq:sigma choice} and \eqref{eq:analyticity constant}, we conclude
\begin{equation}\label{eq:global remainder bound}
\left|\int_0^tR_4(t')\; dt'\right|\le2C_Rc_\ast\delta^{-3/4}E_\ast^2\le\frac14E_\ast.
\end{equation}

Lemma \ref{lem:boundary bound}, \eqref{eq:bootstrap assumption}, and \eqref{eq:sigma boundary smallness} give
\[\begin{aligned}|B(t)|+|B(0)|&\le C_B\sigma\left(E_\sigma(t)^{3/2}+E_\sigma(0)^{3/2}\right)\\&\le C_B\sigma(2^{3/2}+1)E_\ast^{3/2}\\&\le\frac14E_\ast.\end{aligned}\]
Consequently, Proposition \ref{prop:normal form identity} and \eqref{eq:global remainder bound} yield
\[\begin{aligned}E_\sigma(t)&=E_\sigma(0)+B(t)-B(0)+\int_0^tR_4(t')\; dt'\\&\le E_\ast+|B(t)|+|B(0)|+\left|\int_0^tR_4(t')\; dt'\right|\\&\le\frac32E_\ast.\end{aligned}\]
Thus,
\[\sup_{0\le t<T^\ast}E_\sigma(t)\le\frac32E_\ast<2E_\ast.\]
If $T^\ast<T$, applying the uniform local theory at a time in $(\max\{0,T^\ast-\delta/2\},T^\ast)$ extends the solution beyond $T^\ast$, while the continuity of $E_\sigma$ preserves the bound $E_\sigma<2E_\ast$ for a short additional interval. This contradicts the definition of $T^\ast$. Therefore, $T^\ast=T$, and
\[\sup_{0\le t\le T}E_\sigma(t)\le\frac32E_\ast.\]

Applying the same argument to the time-reversed solution gives the corresponding bound on $[-T,0]$. Therefore,
\[u\in C([-T,T];G^{\sigma,0}(\mathbb T)),\]
where $\sigma$ is defined by \eqref{eq:sigma choice}. This proves that, for $s=0$, there exists $c>0$, depending only on $\|u_0\|_{G^{\sigma_0,0}}$ and $\alpha$, such that
\[\sigma(T)=\min\left\{\sigma_0,\frac cT\right\}.\]

\subsection{The general case $s\ge-\alpha/4$}

The reduction to the case $s=0$ is standard, see \cite{SS}. For completeness, we give a short self-contained argument.

Let $u_0\in G^{\sigma_0,s}(\mathbb T)$ and set
\[\rho_0:=\frac{\sigma_0}{2}.\]
A direct computation gives
\[\|u_0\|_{G^{\rho_0,0}}^2 \lesssim_{\sigma_0,s} \|u_0\|_{G^{\sigma_0,s}}^2.\]
Applying the result already proved for $s=0$ with initial radius $\rho_0$, we obtain
\[u\in C([-T,T];G^{\rho(T),0}(\T)),\]
where
\[\rho(T):=\min\left\{\frac{\sigma_0}{2},\frac{c_{\rho_0}}{T}\right\}\]
for some $c_{\rho_0}>0$ depending only on $\|u_0\|_{G^{\sigma_0,s}}$, $\alpha$, $\sigma_0$, and $s$.

For every fixed $\rho>0$, since
\[\|u(t)\|_{G^{\rho/2,s}}^2\lesssim_{\rho,s}\|u(t)\|_{G^{\rho,0}}^2,\]
taking $\rho=\rho(T)$ gives
\begin{equation}\label{eq:general s radius}
u\in C\left([-T,T];G^{\rho(T)/2,s}(\T)\right),
\end{equation}
where
\[\frac{\rho(T)}2=\min\left\{\frac{\sigma_0}{4},\frac{c_{\rho_0}}{2T}\right\}.\]

Let $\delta_s>0$ be the local lifespan provided by Theorem \ref{thm:hks lwp} for $u_0\in G^{\sigma_0,s}(\T)$, and let $u^{(s)}$ denote the corresponding local solution. Set
\[\eta:=\min\left\{\sigma_0,\frac{\rho(\delta_s)}{2}\right\}>0.\]
By \eqref{eq:general s radius},
\[u\in C([-\delta_s,\delta_s];G^{\eta,s}(\T)),\]
while
\[u^{(s)}\in C([-\delta_s,\delta_s];G^{\sigma_0,s}(\T))\hookrightarrow C([-\delta_s,\delta_s];G^{\eta,s}(\T)).\]
Remark \ref{rem:unconditional uniqueness} therefore gives $u=u^{(s)}$ on $[-\delta_s,\delta_s]$.

Set
\[c:=\min\left\{\frac{\sigma_0\delta_s}{4},\frac{c_{\rho_0}}2\right\}.\]
For $0<T\le\delta_s$, we have
\[u\in C([-T,T];G^{\sigma_0,s}(\T)),\]
which implies the desired conclusion. For $T\ge\delta_s$, we have
\[\frac cT\le\min\left\{\frac{\sigma_0}{4},\frac{c_{\rho_0}}{2T}\right\}=\frac{\rho(T)}2,\]
and hence \eqref{eq:general s radius} gives
\[u\in C\left([-T,T];G^{\min\{\sigma_0,c/T\},s}(\T)\right).\]

This completes the proof for solutions with zero spatial mean. The general case follows from Remark \ref{rem:recover mean}.

\begin{rem}\label{rem:recover mean}
The mean-zero assumption causes no loss of generality. Indeed, let
\[m:=\widehat u_0(0),\qquad u_0^\circ:=u_0-m.\]
Then $u_0^\circ$ has zero spatial mean. Let $v$ denote the solution to \eqref{eq:dgkdv} with initial datum $u_0^\circ$, and define
\[u(t,x):=v(t,x-mt)+m.\]
A direct computation shows that $u$ solves \eqref{eq:dgkdv} with initial datum $u_0$. Since the spatial mean is conserved, we have $\widehat v(t,0)=0$. Moreover,
\[\widehat u(t,k)=
\begin{cases}
m,&k=0,\\
e^{-ikmt}\widehat v(t,k),&k\ne0.
\end{cases}\]
Therefore,
\[\begin{aligned}
\|u(t)\|_{G^{\sigma,s}}^2
&=\sum_{k\in\Z}e^{2\sigma|k|}\langle k\rangle^{2s}|\widehat u(t,k)|^2\\
&=|m|^2+\sum_{k\ne0}e^{2\sigma|k|}\langle k\rangle^{2s}|\widehat v(t,k)|^2\\
&=|m|^2+\|v(t)\|_{G^{\sigma,s}}^2.
\end{aligned}\]
Thus, the constant Fourier mode does not affect the radius of spatial analyticity, and the conclusion of Theorem \ref{thm:main} follows from the mean-zero case.
\end{rem}


\begin{thebibliography}{99}

\bibitem{AKS1} J. Ahn, J. Kim, and I. Seo, \emph{On the radius of spatial analyticity for defocusing nonlinear Schr\"odinger equations}, Discrete Contin. Dyn. Syst. 40 (2020), 423--439.

\bibitem{AKS2} J. Ahn, J. Kim, and I. Seo, \emph{Lower bounds on the radius of spatial analyticity for the Kawahara equation}, Anal. Math. Phys. 11 (2021), 28.

\bibitem{AKS3} J. Ahn, J. Kim, and I. Seo, \emph{On the radius of spatial analyticity for the Klein--Gordon--Schr\"odinger system}, J. Differential Equations 321 (2022), 449--474.

\bibitem{AK2026} J. Ahn and C. Kwak, \emph{Persistence of spatial analyticity for the periodic modified KdV equation}, In preparation.

\bibitem{BIT2011} A. V. Babin, A. A. Ilyin, and E. S. Titi, \emph{On the regularization mechanism for the periodic Korteweg--de Vries equation}, Comm. Pure Appl. Math. 64 (2011), no. 5, 591--648.

\bibitem{BC} M. Baldasso and S. Correia, \emph{Nonlinear smoothing implies improved lower bounds on the radius of spatial analyticity for nonlinear dispersive equations}, arXiv preprint arXiv:2507.13083 (2025).

\bibitem{BP} M. Baldasso and M. Panthee, \emph{Improved algebraic lower bound for the radius of spatial analyticity for the generalized KdV equation}, Nonlinear Anal. Real World Appl. 77 (2024), 104054.

\bibitem{Beals} M. Beals, \emph{Self-spreading and strength of singularities for solutions to semilinear wave equations}, Ann. of Math. (2) 118 (1983), no. 1, 187--214.

\bibitem{BG} J. L. Bona and Z. Gruji\'c, \emph{Spatial analyticity properties of nonlinear waves}, Math. Models Methods Appl. Sci. 13 (2003), 345--360.

\bibitem{BGK} J. L. Bona, Z. Gruji\'c, and H. Kalisch, \emph{Algebraic lower bounds for the uniform radius of spatial analyticity for the generalized KdV equation}, Ann. Inst. H. Poincar\'e Anal. Non Lin\'eaire 22 (2005), 783--797.

\bibitem{BGK2} J. L. Bona, Z. Gruji\'c, and H. Kalisch, \emph{Global solutions of the derivative Schr\"odinger equation in a class of functions analytic in a strip}, J. Differential Equations 229 (2006), 186--203.

\bibitem{Bourgain1993-1} J. Bourgain, \emph{Fourier transform restriction phenomena for certain lattice subsets and applications to nonlinear evolution equations. I. Schr\"odinger equations}, Geom. Funct. Anal. 3 (1993), no. 2, 107--156.

\bibitem{Bourgain1993-2} J. Bourgain, \emph{Fourier transform restriction phenomena for certain lattice subsets and applications to nonlinear evolution equations. II. The KdV equation}, Geom. Funct. Anal. 3 (1993), no. 3, 209--262.

\bibitem{Boussinesq} J. Boussinesq, \emph{Th\'eorie des ondes et des remous qui se propagent le long d'un canal rectangulaire horizontal, en communiquant au liquide contenu dans ce canal des vitesses sensiblement pareilles de la surface au fond}, J. Math. Pures Appl. (2) 17 (1872), 55--108.

\bibitem{CKV2022} R. de A. Capistrano--Filho, C. Kwak, and F. J. Vielma Leal, \emph{On the control issues for higher-order nonlinear dispersive equations on the circle}, Nonlinear Anal. Real World Appl. 68 (2022), 103695.

\bibitem{CKSTT} J. Colliander, M. Keel, G. Staffilani, H. Takaoka, and T. Tao, \emph{Sharp global well-posedness for KdV and modified KdV on $\R$ and $\T$}, J. Amer. Math. Soc. 16 (2003), no. 3, 705--749.

\bibitem{Danielski}
A. Danielski, \emph{Analytical Structure of Stationary Flows of an Ideal Incompressible Fluid}, Master's thesis, Concordia University, 2017.

\bibitem{deJager} E. M. de Jager, \emph{On the origin of the Korteweg--de Vries equation}, in A. Festkolloquium ``Rudolf Gorenflo. Fluide aus fraktionaler Sicht''. B. Hans Gebeleins Turbulenz aus stochastischer Sicht, Wellen von Korteweg und de Vries, zellul\"are Diffusion u.a., Berliner Mathematische Gesellschaft, Berlin, 2011, 171--195.

\bibitem{DET} S. Demirbas, M. B. Erdo\u{g}an, and N. Tzirakis, \emph{Existence and uniqueness theory for the fractional Schr\"odinger equation on the torus}, in Some Topics in Harmonic Analysis and Applications, Adv. Lect. Math. (ALM), vol. 34, Int. Press, Somerville, MA, 2016, 145--162.

\bibitem{FH} R. O. Figueira and A. A. Himonas, \emph{Lower bounds on the radius of analyticity for a system of modified KdV equations}, J. Math. Anal. Appl. 497 (2021), 124917.

\bibitem{FNP} R. O. Figueira, M. Nogueira, and M. Panthee, \emph{Lower bounds on the radius of analyticity for a system of nonlinear quadratic interactions of the Schr\"odinger-type equations}, Z. Angew. Math. Phys. 75 (2024), 136.

\bibitem{FP1} R. O. Figueira and M. Panthee, \emph{Decay of the radius of spatial analyticity for the modified KdV equation and the nonlinear Schr\"odinger equation with third order dispersion}, NoDEA Nonlinear Differential Equations Appl. 31 (2024), 68.

\bibitem{FP2} R. O. Figueira and M. Panthee, \emph{Improved lower bound for the radius of analyticity for the modified KdV equation}, arXiv preprint arXiv:2406.08400 (2024).

\bibitem{FP2026} R. O. Figueira and M. Panthee, \emph{Evolution of the radius of analyticity for mKdV-type equations}, Proc. Roy. Soc. Edinburgh Sect. A, First View (2026), 1--27.

\bibitem{GGKM} C. S. Gardner, J. M. Greene, M. D. Kruskal, and R. M. Miura, \emph{Method for solving the Korteweg--de Vries equation}, Phys. Rev. Lett. 19 (1967), 1095--1097.

\bibitem{Get2} T. Getachew, \emph{Asymptotic lower bound on the spatial analyticity radius for solutions of the periodic fifth-order KdV--BBM equation}, Int. J. Differ. Equ. 2025 (2025), 1--9.

\bibitem{Get1} T. Getachew, \emph{Lower bounds of spatial analyticity radius for Benjamin--Bona--Mahony equation on the circle}, J. Math. Model. 13 (2025), 201--208.

\bibitem{BG_NLS} T. Getachew and B. Belayneh, \emph{New asymptotic lower bound for the radius of analyticity of solutions to nonlinear Schr\"odinger equation}, Anal. Appl. (Singap.) 22 (2024), 815--832.

\bibitem{GKO} Z. Guo, S. Kwon, and T. Oh, \emph{Poincar\'e--Dulac normal form reduction for unconditional well-posedness of the periodic cubic NLS}, Comm. Math. Phys. 322 (2013), no. 1, 19--48.

\bibitem{Hadamard1902} J. Hadamard, \emph{Sur les probl\`emes aux d\'eriv\'ees partielles et leur signification physique}, Princeton Univ. Bull. 13 (1902), 49--52.

\bibitem{Hereman} W. Hereman, \emph{Shallow water waves and solitary waves}, in Solitons, Encyclopedia of Complexity and Systems Science Series, M. A. Helal, ed., Springer, New York, 2022, 203--220.

\bibitem{HKS} A. A. Himonas, H. Kalisch, and S. Selberg, \emph{On persistence of spatial analyticity for the dispersion-generalized periodic KdV equation}, Nonlinear Anal. Real World Appl. 38 (2017), 35--48.

\bibitem{HP} A. A. Himonas and G. Petronilho, \emph{Evolution of the radius of spatial analyticity for the periodic BBM equation}, Proc. Amer. Math. Soc. 148 (2020), no. 7, 2953--2967.

\bibitem{Hirayama} H. Hirayama, \emph{Local well-posedness for the periodic higher order KdV type equations}, NoDEA Nonlinear Differential Equations Appl. 19 (2012), no. 6, 677--693.

\bibitem{HK} S. Hong and C. Kwak, \emph{Global well-posedness and nonsqueezing property for the higher-order KdV-type flow}, J. Math. Anal. Appl. 441 (2016), 140--166.

\bibitem{HW} J. Huang and M. Wang, \emph{New lower bounds on the radius of spatial analyticity for the KdV equation}, J. Differential Equations 266 (2019), 5278--5317.

\bibitem{HMW} R. Hunt, B. Muckenhoupt, and R. Wheeden, \emph{Weighted norm inequalities for the conjugate function and Hilbert transform}, Trans. Amer. Math. Soc. 176 (1973), 227--251.

\bibitem{KP} T. Kappeler and J. P\"oschel, \emph{On the periodic KdV equation in weighted Sobolev spaces}, Ann. Inst. H. Poincar\'e Anal. Non Lin\'eaire 26 (2009), 841--853.

\bibitem{KT} T. Kappeler and P. Topalov, \emph{Global wellposedness of KdV in $H^{-1}(\T,\R)$}, Duke Math. J. 135 (2006), no. 2, 327--360.

\bibitem{KatoKawahara} T. Kato, \emph{Low regularity well-posedness for the periodic Kawahara equation}, Differential Integral Equations 25 (2012), no. 11--12, 1011--1036.

\bibitem{KM} T. Kato and K. Masuda, \emph{Nonlinear evolution equations and analyticity. I}, Ann. Inst. H. Poincar\'e Anal. Non Lin\'eaire 3 (1986), 455--467.

\bibitem{K} Y. Katznelson, \emph{An Introduction to Harmonic Analysis}, corrected ed., Dover Publications, Inc., New York, 1976.

\bibitem{KPV1996} C. E. Kenig, G. Ponce, and L. Vega, \emph{A bilinear estimate with applications to the KdV equation}, J. Amer. Math. Soc. 9 (1996), no. 2, 573--603.

\bibitem{KdV} D. J. Korteweg and G. de Vries, \emph{On the change of form of long waves advancing in a rectangular canal, and on a new type of long stationary waves}, Philos. Mag. (5) 39 (1895), no. 240, 422--443.

\bibitem{Kwak2018} C. Kwak, \emph{Periodic fourth-order cubic NLS: Local well-posedness and non-squeezing property}, J. Math. Anal. Appl. 461 (2018), 1327--1364.

\bibitem{KwakKawahara} C. Kwak, \emph{Well-posedness issues on the periodic modified Kawahara equation}, Ann. Inst. H. Poincar\'e Anal. Non Lin\'eaire 37 (2020), no. 2, 373--416.

\bibitem{KY} C. Kwak and C. Yang, \emph{Periodic FPU system: Continuum limit to KdV via regularization and Fourier analysis}, arXiv preprint arXiv:2502.00786 (2025).

\bibitem{KO} S. Kwon and T. Oh, \emph{On unconditional well-posedness of modified KdV}, Int. Math. Res. Not. IMRN 2012 (2012), no. 15, 3509--3534.

\bibitem{Lax} P. D. Lax, \emph{Integrals of nonlinear equations of evolution and solitary waves}, Comm. Pure Appl. Math. 21 (1968), 467--490.

\bibitem{LiShi} J. Li and S. Shi, \emph{Local well-posedness for the dispersion generalized periodic KdV equation}, J. Math. Anal. Appl. 379 (2011), no. 2, 706--718.

\bibitem{LW} K. Liu and M. Wang, \emph{Fixed analytic radius lower bound for the dissipative KdV equation on the real line}, NoDEA Nonlinear Differential Equations Appl. 29 (2022), 57.

\bibitem{Maspero} A. Maspero, \emph{Tame majorant analyticity for the Birkhoff map of the defocusing nonlinear Schr\"odinger equation on the circle}, Nonlinearity 31 (2018), no. 5, 1981--2030.

\bibitem{Miles} J. W. Miles, \emph{The Korteweg--de Vries equation: A historical essay}, J. Fluid Mech. 106 (1981), 131--147.

\bibitem{PS} G. Petronilho and P. L. da Silva, \emph{On the radius of spatial analyticity for the modified Kawahara equation on the line}, Math. Nachr. 292 (2019), 2032--2047.

\bibitem{Rauch}
J. Rauch, \emph{Fourier Series, Integrals, and Sampling From Basic Complex Analysis}, lecture notes, University of Michigan.

\bibitem{Russell} J. S. Russell, \emph{Report on waves}, Report of the Fourteenth Meeting of the British Association for the Advancement of Science, John Murray, London, 1845, 311--390.

\bibitem{Sato} T. Sato, \emph{Spatial analyticity of solutions to quadratic nonlinear Schr\"odinger equations with mass resonance}, J. Math. Anal. Appl. 551 (2025), 129684.

\bibitem{S} S. Selberg, \emph{On the radius of spatial analyticity for solutions of the Dirac--Klein--Gordon equations in two space dimensions}, Ann. Inst. H. Poincar\'e Anal. Non Lin\'eaire 36 (2019), 1311--1330.

\bibitem{SS} S. Selberg and D. O. da Silva, \emph{Lower bounds on the radius of spatial analyticity for the KdV equation}, Ann. Henri Poincar\'e 18 (2017), 1009--1023.

\bibitem{ST} S. Selberg and A. Tesfahun, \emph{On the radius of spatial analyticity for the 1d Dirac--Klein--Gordon equations}, J. Differential Equations 259 (2015), 4732--4744.

\bibitem{ST_gKdV} S. Selberg and A. Tesfahun, \emph{On the radius of spatial analyticity for the quartic generalized KdV equation}, Ann. Henri Poincar\'e 18 (2017), 3553--3564.

\bibitem{Shatah1985} J. Shatah, \emph{Normal forms and quadratic nonlinear Klein--Gordon equations}, Comm. Pure Appl. Math. 38 (1985), 685--696.

\bibitem{TakaokaTsutsumi2004} H. Takaoka and Y. Tsutsumi, \emph{Well-posedness of the Cauchy problem for the modified KdV equation with periodic boundary condition}, Int. Math. Res. Not. IMRN 2004 (2004), 3009--3040.

\bibitem{Tao2001} T. Tao, \emph{Multilinear weighted convolution of $L^2$ functions, and applications to nonlinear dispersive equations}, Amer. J. Math. 123 (2001), no. 5, 839--908.

\bibitem{Tao2006} T. Tao, \emph{Nonlinear Dispersive Equations: Local and Global Analysis}, CBMS Reg. Conf. Ser. Math., vol. 106, American Mathematical Society, Providence, RI, 2006.

\bibitem{TTB} E. Tegegn, A. Tesfahun, and B. Belayneh, \emph{Lower bounds on the radius of spatial analyticity of solution for KdV--BBM type equations}, NoDEA Nonlinear Differential Equations Appl. 30 (2023), 47.

\bibitem{Te} A. Tesfahun, \emph{On the radius of spatial analyticity for cubic nonlinear Schr\"odinger equations}, J. Differential Equations 263 (2017), 7496--7512.

\bibitem{Te2} A. Tesfahun, \emph{Asymptotic lower bound for the radius of spatial analyticity to solutions of KdV equation}, Commun. Contemp. Math. 21 (2019), 1850061.

\bibitem{W} M. Wang, \emph{Nondecreasing analytic radius for the KdV equation with a weakly damping}, Nonlinear Anal. 215 (2022), 112653.

\bibitem{ZF} V. E. Zakharov and L. D. Faddeev, \emph{Korteweg--de Vries equation: A completely integrable Hamiltonian system}, Funct. Anal. Appl. 5 (1971), no. 4, 280--287.

\bibitem{ZLD} Z. Zhang, Z. Liu, and Y. Deng, \emph{Lower bounds on the radius of spatial analyticity for the higher order nonlinear dispersive equation on the real line}, Discrete Contin. Dyn. Syst. Ser. B 29 (2024), 937--970.

\end{thebibliography}
\end{document}